\documentclass[pdflatex, sn-mathphys]{sn-jnl}

\jyear{2022}%

\theoremstyle{thmstyleone}%
\newtheorem{theorem}{Theorem}
\newtheorem{proposition}[theorem]{Proposition}%
\newtheorem{corollary}[theorem]{Corollary}
\newtheorem{lemma}[theorem]{Lemma}

\usepackage{xcolor}
\definecolor{VivaMagenta}{HTML}{BE3455}

\theoremstyle{thmstyletwo}%
\newtheorem{example}{Example}%
\newtheorem{remark}{Remark}%
\newtheorem{question}{Question}

\theoremstyle{thmstylethree}%
\newtheorem{definition}{Definition}%

\usepackage{lmodern}
\usepackage{mathtools}
\usepackage[caption=false]{subfig}
\catcode`\|=12\relax
\usepackage{tikz}

\newcommand{\bbC}{\mathbb{C}}

\renewcommand{\R}{\mathbb{R}}
\renewcommand{\Z}{\mathbb{Z}}

\begin{document}

\title[Uniqueness from samples and stability in Gabor phase retrieval]{On the connection between uniqueness from samples and stability in Gabor phase retrieval}


\author[1]{\fnm{Rima} \sur{Alaifari}}\email{rima.alaifari@math.ethz.ch}

\author*[2]{\fnm{Francesca} \sur{Bartolucci}}\email{F.bartolucci@tudelft.nl}

\author[3]{\fnm{Stefan} \sur{Steinerberger}}\email{steinerb@uw.edu}

\author[4]{\fnm{Matthias} \sur{Wellershoff}}\email{wellersm@umd.edu}

\affil[1]{\orgname{ETH Z\"urich}, \orgdiv{Department of Mathematics, Seminar for Applied Mathematics}, \orgaddress{\street{R\"amistrasse 101}, \postcode{CH-8092} \city{Z\"urich}, \country{Switzerland}}}

\affil[2]{\orgname{TU Delft}, \orgdiv{Institute of Applied Mathematics}, \orgaddress{\street{Mekelweg 4}, \postcode{2628 CD} \city{Delft}, \country{The Netherlands}}}

\affil[3]{\orgname{University of Washington}, \orgdiv{Department of Mathematics}, \orgaddress{\street{C-138 Padelford}, \city{Seattle}, \state{WA} \postcode{98195-4350}, \country{USA}}}

\affil[4]{\orgname{University of Maryland}, \orgdiv{Department of Mathematics}, \orgaddress{\street{4176 Campus Drive}, \city{College Park}, \state{MD} \postcode{20742}, \country{USA}}}


\abstract{
\emph{Gabor phase retrieval} is the problem of reconstructing a signal from only the magnitudes of its Gabor transform. Previous findings suggest a possible link between unique solvability of the discrete problem (recovery from measurements on a lattice) and stability of the continuous problem (recovery from measurements on an open subset of $\mathbb{R}^2$). In this paper, we close this gap by proving that such a link cannot be made. More precisely, we establish the existence of functions which break uniqueness from samples without affecting stability of the continuous problem. Furthermore, we prove the novel result that counterexamples to unique recovery from samples are \emph{dense} in $L^2(\mathbb{R})$.
Finally, we develop an intuitive argument on the connection between directions of instability in phase retrieval and certain Laplacian eigenfunctions associated to small eigenvalues.
}




\maketitle

\section{Introduction}\label{sec:intro}

\emph{Phase retrieval} is a broad term encompassing many different inverse problems in imaging and signal processing, in which one seeks to reconstruct an object from measurements that arise as magnitudes of some quantities and lack phase information. While in imaging applications such as coherent diffraction imaging, this loss of phase information is due to the physics of the data acquisition process, in audio processing applications (cf.~\cite{prusa2017phase}) working with phaseless information is often a choice (to avoid recovery from noisy phases for example). 

A particular instance of such a phase recovery problem is \emph{STFT phase retrieval,} i.e. the recovery of a signal $f \in L^2(\R)$ from the magnitude of its short-time Fourier transform with respect to a window $\psi \in L^2(\R)$, which is defined as
\begin{equation}\label{eqn:STFT}
    \mathcal{V}_\psi f(x,\omega) := \int_{\R} f(t) \overline{\psi(t-x)} \mathrm{e}^{-2\pi \mathrm{i} t \omega} \,\mathrm{d} t, \qquad (x,\omega) \in \R^2.
\end{equation}
STFT phase retrieval then refers to the problem of recovering $f$ from $\left| V_\psi f\right(x,\omega)|,$ up to a constant global phase factor (cf.~equation~\eqref{eq:uptoglobalphase}).

A popular window choice for the STFT is the Gauss function $\varphi(t)=2^{1/4} \mathrm{e}^{-\pi t^2}$, in which case one also refers to the STFT as the \emph{Gabor transform} defined as
\begin{equation}\label{eqn:Gabor}
    \mathcal{G} f(x,\omega) := 2^{1/4} \int_{\R} f(t) \mathrm{e}^{-\pi(t-x)^2} \mathrm{e}^{-2\pi \mathrm{i} t \omega} \,\mathrm{d} t, \qquad (x,\omega) \in \R^2,
\end{equation}
and the related STFT phase retrieval as \emph{Gabor phase retrieval}.

In what follows, we will distinguish two setups: first, \emph{(continuous)} Gabor phase retrieval problems, in which one seeks to recover $f$ from $\mathcal{A}_\Omega(f):=(\lvert \mathcal{G} f (x,\omega)\rvert)_{(x,\omega) \in \Omega}$, for $\Omega$ containing an open subset of $\R^2$. Second, and more relevant in view of applications, \emph{sampled} Gabor phase retrieval problems that deal with the recovery of $f$ (up to global phase) from $\mathcal{A}_\Lambda(f) := (\lvert \mathcal{G} f (x,\omega)\rvert)_{(x,\omega) \in \Lambda}$, where $\Lambda$ is a discrete subset of the time-frequency plane $\R^2$.

It is by now well-known that $f$ may be uniquely recovered (up to global phase) from $\mathcal{A}_{\R^2}(f)$ and that this recovery is weakly but not strongly stable, i.e.~the inverse phase retrieval operator, $\mathcal{A}_{\R^2}^{-1}$, is continuous but not uniformly continuous \cite{alaifari2017phase}. In view of this insight, one can attempt to derive upper bounds for the local Lipschitz constants $c_{\R^2}(f)$ for the inverse phase retrieval operator (cf.~equation \eqref{eq:local_stability_constant}). 

\begin{figure}
    \begin{center}
            \includegraphics[width=.49\textwidth]{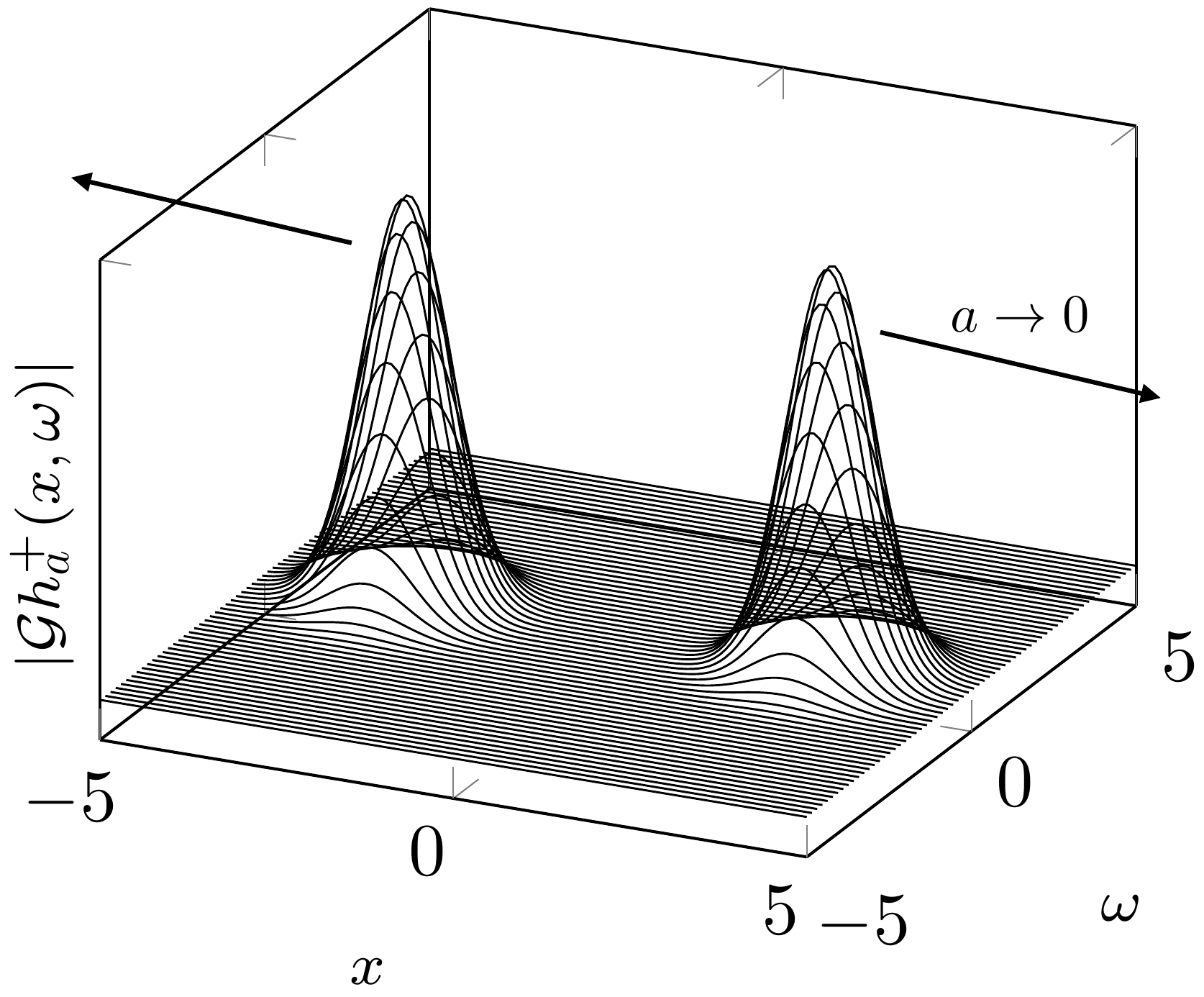}
        \label{fig:introduction}
    \caption{The Gabor transform magnitudes $|\mathcal{G}h_a^+|$ for a=1/6. The two bumps move apart as the sample rate $a$ goes to zero. As a consequence, the local Lipschitz constant of $h_a^+$ increases as $a$ goes to zero.}    
    \end{center}
\end{figure}

Different further aspects of uniqueness and stability for Gabor phase retrieval have been studied over the past decade \cite{iwen2023phase,eldar2014sparse,bojarovska2016phase,li2017phase,pfander2019robust,salanevich2023injectivity,alaifari2019stable,alaifari2021Gabor,alaifari2021stability,alaifari2021phase,alaifari2021uniqueness,grohs2023injectivity,grohs2019stable,grohs2021stable,grohs2022foundational,grohs2022non,grohs2022phaseless,grohs2022multi,wellershoff2023sampling,wellershoff2022phase,wellershoff2021injectivity}. Here, we want to highlight two findings that motivate the study of this paper:
\begin{itemize}
    \item Prior work by two authors of this paper on the (non-)uniqueness of sampled Gabor phase retrieval \cite{alaifari2021phase} shows that sampled Gabor phase retrieval does not enjoy uniqueness (for signals in $L^2(\R)$) when the sampling set $\Lambda$ is any (shifted) lattice in $\R^2$. This is done by constructing explicit counterexamples:

    \begin{theorem}[{\cite[Theorem 1 on p. 6]{alaifari2021phase}}]\label{thm:original_counterexamples}
        For any $a>0$, the functions 
        \begin{equation*}
            h^\pm_a(t) := 2^{1/4} \exp(-\pi t^2) \left( \cosh\left( \frac{\pi t}{a} \right) \pm \mathrm{i} \sinh\left( \frac{\pi t}{a} \right) \right), \qquad t \in \R,
        \end{equation*}
         do \emph{not} agree (up to a constant global phase factor) and yet
        \begin{equation*}
            |\mathcal{G}h_a^+|=|\mathcal{G}h_a^{-}|\qquad \text{on}\quad \R\times a\mathbb{Z}\supset a\mathbb{Z}^2.
        \end{equation*}
    \end{theorem}

    As a consequence, a priori knowledge about the signals is necessary to restore uniqueness in sampled Gabor phase retrieval, and the search for proper subspaces of $L^2(\R)$ enjoying uniqueness from sampled Gabor transform magnitudes has attracted recent attention  \cite{alaifari2021uniqueness, grohs2023injectivity,wellershoff2021injectivity}.

    \item When the sampling lattice is sufficiently fine (i.e. $a\to 0$), the counterexamples $h^\pm_a$ strongly resemble signals proposed by Grohs and one of the authors to demonstrate that continuous Gabor phase retrieval is severely ill-posed \cite{alaifari2021Gabor}: The magnitude measurements $|\mathcal{G}h_a^\pm|$ are concentrated on disjoint domains in the time-frequency plane and this results in a source of instability for Gabor phase retrieval (i.e. the local Lipschitz constants $c_{\R^2}(h^\pm_a)$ are large). More precisely, the local stability constant of the counterexamples increases as the sampling lattice becomes finer and finer (see Figure~\ref{fig:introduction}).
\end{itemize}

We therefore observe that function pairs breaking stability of continuous Gabor phase retrieval seemingly resemble function pairs breaking uniqueness of sampled Gabor phase retrieval. One is thus tempted to expect a direct connection between uniqueness of sampled Gabor phase retrieval and stability of continuous Gabor phase retrieval. This observation naturally leads to the following central question:

\begin{question}\label{conjecture}
    Consider the set $\mathcal{M}_{\nu}(\R^2)$ of signals whose local Lipschitz constant (defined in equation~\eqref{eq:local_stability_constant2}) is upper bounded by $\nu > 0$. Is there a lattice $\Lambda \subset \R^2$ such that $\lvert \mathcal{G} f \rvert = \lvert \mathcal{G} g \rvert$ on $\Lambda$ implies that there exists an $\alpha \in \R$ with $f = \mathrm{e}^{\mathrm{i} \alpha} g$ for all $f,g \in \mathcal{M}_{\nu}(\R^2)$?
\end{question}

    

\begin{figure}
    \begin{center}
        \includegraphics[width=.49\textwidth]{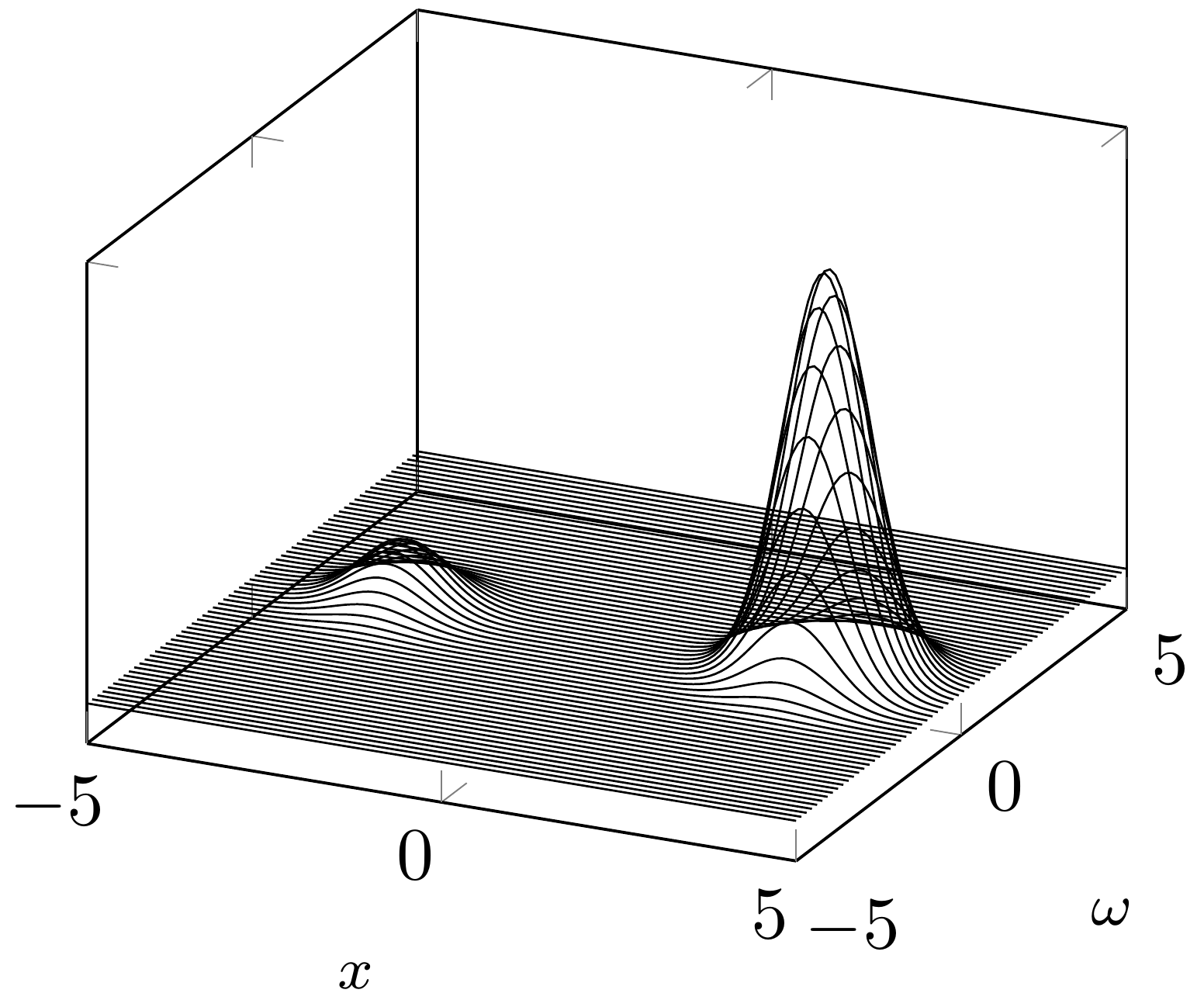}
        \caption{A plot of the Gabor transform magnitude of the counterexamples whose existence is postulated in Theorem~\ref{thm:main}.}   
        \label{fig:introduction2}
    \end{center}
\end{figure}

This question has been the starting point and main motivation of this paper. Along the way, we have established other results that are interesting on their own. The four main contributions of our work can be summarized as follows:
\begin{enumerate}
    \item \emph{The set of counterexamples is dense in $L^2(\R)$.} For every lattice and for every square-integrable function $f$, we can construct signals which are arbitrarily close to $f$, do not agree up to a global constant phase factor and yet have Gabor transform magnitudes agreeing on the lattice (see Section \ref{sec:classcounterexamplesdense}). 

    \item \emph{The notion of uniqueness for sampled Gabor phase retrieval is fragile.} The Gaussian can be recovered from sampled Gabor transform magnitudes on sufficiently fine lattices but 
    for each of such lattices there exist counterexamples which are arbitrarily close to the Gaussian (cf. Section \ref{sec:fragility}).
    
    \item \emph{Considering the class of signals for which the local Lipschitz constant satisfies a uniform bound is not enough to restore uniqueness in sampled Gabor phase retrieval.} Much to our surprise, we answer Question \ref{conjecture} in the negative. This is due to the existence of counterexamples which break uniqueness without affecting stability, see Figure~\ref{fig:introduction2}. Indeed, the counterexamples in Theorem \ref{thm:main} are constructed to be arbitrarily close to the normalized Gaussian $\varphi(t) = 2^{1/4} \mathrm{e}^{-\pi t^2}$ which is known to enjoy very strong stability properties. Restoring uniqueness will necessitate a more restrictive signal class (cf. Section \ref{sec:stability}). 
    
    \item \emph{Local stability of (Gabor) phase retrieval may be quantified by Laplacian eigenvalues.} In particular, we suggest that small Laplacian eigenvalues correspond to unstable directions: If there are only very few small eigenvalues of the Laplacian, there are only few directions of instability. Moreover, each direction of instability corresponds to an associated Laplacian eigenfunction (see Section \ref{sec:further_remarks}).
\end{enumerate}

\begin{remark}
    As there is a lot of notation introduced throughout this paper, we include a list of symbols in Appendix~\ref{app:A} for the convenience of the reader.
\end{remark}

\section{Preliminaries and main concepts}\label{sec:preliminaries}

We will consider signals in the \emph{modulation spaces}
\begin{equation}\label{eq:modspace}
    M^p(\R) := \left\{ f \in \mathcal{S}'(\R) : \mathcal{G} f \in L^p(\R^2) \right\}, \qquad 1 \leq p \leq \infty,
\end{equation}
where $\mathcal{S}'(\R)$ denotes the class of tempered distributions. $M^p(\R)$ can be equipped with the norm
\begin{equation*}
    \lVert f \rVert_{M^p(\R)} := \lVert \mathcal{G} f \rVert_p.
\end{equation*}
Let us emphasize that we are exclusively interested in the modulation spaces with parameter $p \in [1,2]$. We will thus always be working with functions rather than  abstract distributions. In particular, the following simple inclusion holds (cf.~\cite[Proposition~1.7 on p.~408]{toft2004continuity}):
\begin{equation*}
    M^p(\R) \subset L^r(\R), \qquad r \in [p,p'],
\end{equation*}
where $p' \in [2,\infty]$ denotes the H{\"o}lder conjugate of $p$ and where we have equality as sets if $p = 2$. We observe that the above inclusion implies that $M^p(\R) \subset L^2(\R)$, for $p \in [1,2]$, such that the application of the Gabor transform to signals $f \in M^p(\R)$ is well-defined.

Next, we note that one cannot distinguish between $f$ and $\mathrm{e}^{\mathrm{i} \alpha} f$ on the basis of the magnitudes of the Gabor transform. As is usual in phase retrieval, we will therefore only seek to recover $f$ \emph{up to a global constant phase factor}. To formalize this, we introduce the equivalence relation 
\begin{equation}\label{eq:uptoglobalphase}
    f \sim g :\iff \exists \, \alpha \in \R: f = \mathrm{e}^{\mathrm{i} \alpha} g
\end{equation}
on $M^p(\R)$.

For ease of notation, we introduce the following operations in the context of $L^p(\mathbb{R})$ with $p \in [1,\infty]$: The translation operator $\operatorname{T}_x : L^p(\mathbb{R}) \to L^p(\mathbb{R})$ is defined as
\begin{equation}\label{eq:definitionoftranslation}
    \operatorname{T}_x f(t) := f(t-x), \quad t \in \mathbb{R},
\end{equation}
for $x \in \mathbb{R}$; the modulation operator $\operatorname{M}_\omega : L^p(\mathbb{R}) \to L^p(\mathbb{R})$ is defined as
\begin{equation*}
    \operatorname{M}_\omega f(t) := f(t) \mathrm{e}^{2\pi\mathrm{i} t \omega}, \quad t \in \mathbb{R},
\end{equation*}
for $\omega \in \mathbb{R}$.

Moreover, to represent rotation in $\mathbb{R}^2$ by an angle $\theta \in \mathbb{R}$, we use the rotation operator $\operatorname{R}_\theta : \mathbb{R}^2 \to \mathbb{R}^2$, which, in matrix form, can be expressed as
\begin{equation}\label{eq:definitionofrotation}
    \operatorname{R}_\theta := \begin{pmatrix}
    \cos \theta & - \sin \theta \\
    \sin \theta & \cos \theta
    \end{pmatrix}.
\end{equation}

Finally, we say that a \emph{lattice} $\Lambda \subset \mathbb{R}^2$ is a discrete subset of the time-frequency plane that can be expressed as $L \mathbb{Z}^{k}$, where $L \in \mathbb{R}^{2 \times k}$ is a matrix with linearly independent columns and $k \in \{1,2\}$.

\subsection{The Fock space and the Bargmann transform}\label{ssec:fockspace}

Throughout this work, and notably in Section~\ref{sec:classcounterexamplesdense}, we rely on the well-known connection between the Gabor transform, the Bargmann transform and the Fock space \cite[Section 3.4 on pp.~53--58]{groechenig2001foundations}. The \emph{Fock space} $\mathcal{F}^2(\bbC)$ is the Hilbert space of all entire functions for which the norm induced by the inner product
\begin{equation*}
    (F,G)_\mathcal{F} := \int_\bbC F(z) \overline{G(z)} \mathrm{e}^{-\pi \lvert z \rvert^2} \,\mathrm{d} z
\end{equation*}
is finite. Interestingly, the Fock space is isomorphic to $L^2(\R)$ with the \emph{Bargmann transform} $\mathcal{B} : L^2(\R) \to \mathcal{F}^2(\bbC)$, given by
\begin{equation*}
    \mathcal{B} f (z) := 2^{1/4} \int_\R f(t) \mathrm{e}^{2 \pi t z - \pi t^2 - \frac{\pi}{2} z^2} \,\mathrm{d} t, \qquad z \in \bbC,
\end{equation*}
acting as the isomorphism. Moreover, the Bargmann transform and the Gabor transform are intimately related via the formula
\begin{equation}\label{eq:Gabor_Bargmann}
    \mathcal{G} f(x,-\omega) = \mathrm{e}^{\pi\mathrm{i}x\omega} \mathcal{B} f(x + \mathrm{i} \omega) \mathrm{e}^{-\frac{\pi}{2}(x^2 + \omega^2)}, \qquad (x,\omega) \in \R^2.
\end{equation}
Therefore, the range of the Gabor transform can be identified with $\mathcal{F}^2(\bbC)$ by reflection and multiplication with a smooth, non-zero function.

For our purposes, it is interesting to note that the Fock space is the space of all entire functions with a specific growth. Let us briefly explain this: The \emph{order} $\rho \geq 0$ of an entire function $F$ is the infimum of all $r > 0$ for which $F(z) \in \mathcal{O}(\exp(\lvert z \rvert^r))$ as $\lvert z \rvert \to \infty$. In the case where $F$ is an entire function of order $\rho \in (0,\infty)$, the \emph{type} $\tau \geq 0$ of $F$ is the infimum of all $t > 0$ such that $F(z) \in \mathcal{O}(\exp(t \lvert z \rvert^\rho))$ as $\lvert z \rvert \to \infty$. An entire function $F$ with order $\rho < 1$ or order $\rho = 1$ and type $\tau < \infty$ is said to be of \emph{exponential type}.

According to \cite[Theorem~3.4.2 on p.~54]{groechenig2001foundations}, all functions in the Fock space are either of order $\rho < 2$ or of order $\rho = 2$ and type $\tau \leq \pi/2$. Conversely, it is readily seen that functions of order $\rho < 2$ and functions of order $\rho = 2$ with type $\tau < \pi/2$ belong to the Fock space. In particular, functions of exponential type are in the Fock space. However, functions of order $\rho = 2$ and type $\tau = \pi/2$ must not necessarily belong to the Fock space as demonstrated in \cite{beneteau2010extremal}.

\subsection{Laplacian eigenvalues, the Poincar{\'e}, and the Cheeger constant}\label{ssec:prelim_LaplacePoincareCheeger}

In their work, Grohs and Rathmair \cite{grohs2019stable} establish a connection between the Cheeger constant from spectral geometry and the local stability of the Gabor phase retrieval problem. They do so by linking stability to the Poincar{\'e} constant, which is the reciprocal of the first non-trivial eigenvalue of the Laplace operator. The connection to the Cheeger constant is then made through Cheeger's inequality (cf.~equation~\eqref{eq:Cheeger}). Since we use the link between the local stability of the Gabor phase retrieval problem and the Poincar{\'e} constant in Section~\ref{sec:stability}, and refine the findings in \cite{grohs2019stable} by making use of the Laplacian eigenvalues in Section~\ref{sec:further_remarks}, we briefly explain these concepts in the following.

\subsubsection*{Laplacian eigenvalues, the Poincar{\'e}, and the Cheeger constant in Riemannian geometry}

To facilitate understanding, let us begin by briefly summarizing some fundamental concepts from spectral geometry \cite{chavel1984eigenvalues}. For a compact, connected smooth Riemannian manifold $M$ of dimension $n$, with Laplace--Beltrami operator $\Delta = -\operatorname{div} \operatorname{grad}$, we consider the following eigenvalue problem: Find all $\lambda \in \mathbb{R}$ such that there exists a non-trivial solution $\phi \in C^2(M)$ to $\Delta \phi = \lambda \phi$. It is known that the eigenvalues form a sequence $0 = \lambda_0 < \lambda_1 < \dots$, where the non-negativity of the eigenvalues follows from the divergence theorem which implies\footnote{Here, $\mathrm{d} V$ denotes integration with respect to the Riemannian measure.}
\begin{equation*}
    \lambda_i \cdot \int_M \lvert \phi_i \rvert^2 \,\mathrm{d} V = \int_M \lvert \operatorname{grad} \phi_i \rvert^2 \,\mathrm{d} V,
\end{equation*}
for all eigenpairs $(\lambda_i,\phi_i)$. This equation also demonstrates that the only eigenfunctions with trivial eigenvalue are the constant functions.

Let $H^1(M)$ denote the Sobolev space of $L^2$-functions with square-integrable weak derivatives up to first order. It can be shown that
\begin{equation}\label{eq:Poincare_manifold}
    \int_M \lvert f \rvert^2 \,\mathrm{d} V \leq \lambda_1^{-1} \cdot \int_M \lvert \operatorname{grad} f \rvert^2 \,\mathrm{d} V
\end{equation}
for all $f \in H^1(M)$ with zero mean. We refer to this inequality as the \emph{Poincar{\'e} inequality} and to the reciprocal $\lambda_1^{-1/2}$ as the \emph{Poincar{\'e} constant}. It is worth noting that equality in equation~\eqref{eq:Poincare_manifold} is achieved by the eigenfunctions with eigenvalue $\lambda_1$, making $\lambda_1^{-1}$ the smallest constant for which equation~\eqref{eq:Poincare_manifold} holds.

Lastly, the first non-trivial eigenvalue of the Laplace operator is connected to \emph{Cheeger's (isoperimetric) constant}\footnote{Here, $V(\cdot)$ denotes $n$-dimensional volume, while $A(\cdot)$ represents $(n-1)$-dimensional volume.},
\begin{equation*}
    h(M) := \inf_S \frac{A(S)}{\min\{ V(M_1), V(M_2)\}},
\end{equation*}
where the infimum is taken over all compact $(n-1)$-dimensional submanifolds $S$ of $M$ that divide $M$ into two open submanifolds $M_1$ and $M_2$ satisfying $\partial M_1 = \partial M_2 = S$: Specifically, we have \cite{cheeger1971lower}
\begin{equation}\label{eq:Cheeger}
    \lambda_1 \geq h(M)^2/4.
\end{equation}
Intuitively, the Cheeger constant measures the ease with which we can cut the manifold $M$ into two parts of roughly equal size. The Calabi dumbbell (cf.~Fig.~\ref{fig:dumbbell}) nicely exemplifies this: It has a small Cheeger constant because cutting through the middle (the "bridge") requires only a short incision to separate the two large regions on either side of the bridge.

\subsubsection*{The Poincar{\'e} constant and stability in Gabor phase retrieval}

Next, we provide a brief overview of the connection between the weighted Poincaré constant and stability of the Gabor phase retrieval problem as established in \cite{grohs2019stable}.

Let $1 \leq p < \infty$, $\Omega \subseteq \mathbb{R}^2$ be a domain, and $w$ be a weight on $\Omega$. In the following, $L^p(\Omega,w)$ and $W^{k,p}(\Omega,w)$ denote the respective weighted spaces on $\Omega$. The \emph{weighted Poincaré constant} is defined as
\begin{multline}\label{eq:poincdefnsup}
    C_{\mathrm{poinc}}(p,\Omega,w) := \sup \bigg\{ \frac{\lVert F - F_\Omega^w \rVert_{L^p(\Omega,w)}}{\lVert \nabla F \rVert_{L^p(\Omega,w)}} : \\
    F \in W^{1,p}(\Omega,w) \cap \mathcal{M}(\Omega),~F \neq \mathrm{const.} \bigg\}
\end{multline}
in \cite{grohs2019stable}, where
\begin{equation*}
    F_\Omega^w := \frac{1}{w(\Omega)}\int_{\Omega} F(x)w(x) \mathrm{d} x, \qquad w(\Omega) := \int_{\Omega} w(x) \mathrm{d} x,
\end{equation*}
and $\mathcal{M}(\Omega)$ denotes the set of functions $F: \Omega \to \mathbb{C}$ for which $\{x+iy : (x,y) \in \Omega\} \ni (x+iy) \mapsto F(x,y) \in \mathbb{C}$ is meromorphic.

\begin{remark}
    It is worth noting that the definition of the Poincaré constant used here slightly deviates from the classical definition, as we specifically consider the case for which results from \cite{grohs2019stable} are applicable.
\end{remark}

The central result of \cite{grohs2019stable} states that under certain assumptions (see Theorem 5.9 therein), the following stability result holds: For every function $f \in M^{p}(\mathbb{R})$, there exists a constant $c > 0$ such that for every $g \in M^{p}(\mathbb{R})$,
\begin{equation}\label{eq:stab-GR}
    \inf_{\alpha \in \mathbb{R}} \lVert \mathcal{G}f-\mathrm{e}^{\mathrm{i}\alpha}\mathcal{G}g\rVert_{L^p(\Omega)} \leq c\ (1+C_{\mathrm{poinc}}(p,\Omega,\lvert \mathcal{G}f \rvert^p)) \cdot \lVert \lvert \mathcal{G}f \rvert - \lvert \mathcal{G}g \rvert \rVert_{\mathcal{D}_{p,q}^{1,4}(\Omega)},
\end{equation}
where
\begin{equation}\label{eq:normmeasuremnts}
    \lVert F \rVert_{\mathcal{D}_{p,q}^{1,4}(\Omega)} := \lVert F \rVert_{W^{1,p}(\Omega)} + \lVert F \rVert_{L^q(\Omega)}
    + \lVert (x,\omega) \mapsto (|x|+|\omega|)^4F(x,\omega)\rVert_{L^q(\Omega)}.
\end{equation}

\subsubsection*{The Cheeger constant and stability in Gabor phase retrieval}

Finally, the Poincar{\'e} constant can be upper bounded in terms of the \emph{Cheeger constant},
\begin{equation*}
    h_{p,\Omega}(f) := \inf_S \frac{\lVert \mathcal{G} f \rVert_{L^p(\partial S)}^p}{\min\{ \lVert \mathcal{G} f \rVert_{L^p(S)}^p, \lVert \mathcal{G} f \rVert_{L^p(S^\mathrm{c})}^p \}},
\end{equation*}
where the infimum is taken over all open subsets $S$ of $\Omega$ for which $\partial S \cap \Omega$ is smooth \cite{grohs2019stable}. Here, $f \in M^p(\mathbb{R})$ for some $p \in [1,2]$ and $\Omega \subset \mathbb{R}^2$ is a connected domain. The Cheeger constant quantifies the ease with which we can cut $|\mathcal{G} f|$ into two parts of roughly equal size. Consider the function $h_a^+$ whose Gabor transform magnitude is depicted in Figure~\ref{fig:introduction}, for instance: Cutting right down the middle of the two bumps is easy, i.e., $\lVert \mathcal{G} f \rVert_{L^p(\partial S)}$ is small. It can then be shown that the Cheeger constant of $h_a^+$ is small. We can compare this to the function whose Gabor transform magnitude is depicted in Figure~\ref{fig:introduction2}: There, cutting right down the middle of the two bumps is even easier. However, when one of the bumps is significantly smaller than the other, $\min\{ \lVert \mathcal{G} f \rVert_{L^p(S)}, \lVert \mathcal{G} f \rVert_{L^p(S^\mathrm{c})} \}$ becomes small as well, and it is not clear whether $h_{p,\Omega}(f)$ is small.


\section{Uniqueness is fragile in sampled Gabor phase retrieval}\label{sec:fragility}

In this section, we show that there exist functions which can be uniquely recovered from sampled Gabor transform magnitudes but that, for every function $f \in L^2(\R)$, there exist counterexamples to sampled Gabor phase retrieval which are arbitrarily close to $f$. Therefore, \emph{uniqueness is a fragile concept in sampled Gabor phase retrieval}.

\subsection{On counterexamples for sampled Gabor phase retrieval}
\label{ssec:counterexamples}

Let us first define what is meant by counterexamples for sampled Gabor phase retrieval.

\begin{definition}\label{def:counterexamples}
    Let $\Lambda \subset \R^2$. The class of \emph{counterexamples for sampled Gabor phase retrieval} on $\Lambda$ is defined by
    \begin{equation*}
        \mathfrak{C}(\Lambda) := \left\{ f \in L^2(\R) \,:\, \exists \, g \in L^2(\R) \mbox{ s.t.~} f\not\sim g \mbox{ and } \lvert \mathcal{G} f \rvert = \lvert \mathcal{G} g \rvert \mbox{ on } \Lambda \right\}.
    \end{equation*}
    An element $f \in \mathfrak{C}(\Lambda)$ is called a \emph{counterexample} for sampled Gabor phase retrieval on $\Lambda$.
\end{definition} 

There is a principled way of generating large families of counterexamples starting with two counterexamples only. Specifically, we may consider $f,g \in L^2(\R)$ such that $f$ and $g$ do not agree up to global phase while 
\begin{equation}\label{eq:magnitudes_agree}
    \lvert \mathcal{G} f \rvert = \lvert \mathcal{G} g \rvert \mbox{ on } \Lambda,
\end{equation}
where $\Lambda \subset \R^2$. If for some element $m \in \mathcal{F}^2(\mathbb{C})$ of the Fock space, it is true that $m \cdot \mathcal{B} f, m \cdot \mathcal{B} g \in \mathcal{F}^2(\bbC)$, i.e.~both $m \cdot \mathcal{B} f$ and $m \cdot \mathcal{B} g$ are in the Fock space, then $f_m := \mathcal{B}^{-1} (m \cdot \mathcal{B} f)$ and $g_m := \mathcal{B}^{-1} (m \cdot \mathcal{B} g)$ are well-defined functions in $L^2(\R)$ such that $\lvert \mathcal{G} f_m \rvert = \lvert \mathcal{G} g_m \rvert$ on $\Lambda$. Additionally, $f_m$ and $g_m$ will generally not agree up to global phase. The above construction works because the magnitudes of the Gabor transforms of $f_m$ and $g_m$ agree on $\Lambda$ if and only if the magnitudes of their Bargmann transforms agree. By construction, their Bargmann transforms are $m \cdot \mathcal{B} f$ and $m \cdot \mathcal{B} g$ whose magnitudes must agree on $\Lambda$ by equation~\eqref{eq:magnitudes_agree}. While this idea can be explored in various ways, we choose to focus on its most straightforward corollary in this context, applying it to the counterexamples presented in \cite{alaifari2021phase}.


Let $a > 0$ and denote the normalized Gaussian by $\varphi(t) = 2^{1/4} \exp(-\pi t^2)$ for $t \in \R$. Recall from the introduction (cf.~Theorem~\ref{thm:original_counterexamples}) that the functions
\begin{equation}
    \label{eq:original_counterexamples}
    h^\pm(t) := \varphi(t) \left( \cosh\left( \frac{\pi t}{a} \right) \pm \mathrm{i} \sinh\left( \frac{\pi t}{a} \right) \right), \qquad t \in \R,
\end{equation}
do not agree up to global phase and still satisfy $\lvert \mathcal{G} h^+ \rvert = \lvert \mathcal{G} h^- \rvert$ on $\R \times a \mathbb{Z}$.



The Bargmann transforms of the signals $h^\pm \in L^2(\R)$ are given by 
\begin{equation*}
    \mathcal{B} h^\pm (z) = \mathrm{e}^{\tfrac{\pi}{8 a^2}} \left( \cosh \left(\frac{\pi z}{2a}\right) \pm \mathrm{i} \sinh \left(\frac{\pi z}{2a}\right) \right), \qquad z \in \bbC,
\end{equation*}
and are thus of exponential type. It therefore follows that, if $m \in \mathcal{F}^2(\bbC)$ is not of second order and type $\pi/2$\footnote{If $m$ is of second order and type $\pi/2$, then $m \cdot \mathcal{B} h^\pm$ is not guaranteed to be in the Fock space (cf.~Remark~\ref{rem:interesting_examples}).}, then $m \cdot \mathcal{B} h^\pm \in \mathcal{F}^2(\bbC)$ and thus 
\begin{equation*}
    h_m^\pm := \mathcal{B}^{-1} ( m \cdot \mathcal{B} h^\pm ) \in L^2(\R)
\end{equation*}
satisfy $\lvert \mathcal{G} h_m^+ \rvert = \lvert \mathcal{G} h_m^- \rvert$ on $\R \times a \mathbb{Z}$. In particular, we can consider $m_\tau(z) := \exp(\pi \tau z)$, for $\tau \in \R$ and $z \in \bbC$, which is of exponential type and therefore in the Fock space. The corresponding functions 
\begin{equation*}
    h_\tau^\pm := \mathcal{B}^{-1} ( \mathrm{e}^{\pi \tau \cdot} \cdot \mathcal{B} h^\pm ) \in L^2(\R)
\end{equation*}
do not agree up to global phase and have Gabor transform magnitudes that agree on $\R \times a \mathbb{Z}$. The former is true because $m_\tau = \exp(\pi \tau \cdot)$ has no roots while $\mathcal{B} h^+$ and $\mathcal{B} h^-$ have disjoint root sets such that $m_\tau \cdot \mathcal{B} h^+ \not\sim m_\tau \cdot \mathcal{B} h^-$. Therefore, the linearity of the Bargmann transform implies that $h_\tau^+ \not\sim h_\tau^-$. We visualize the Gabor transform magnitude of $h_\tau^+$ in Figure~\ref{fig:introduction2} and note that one of the two bumps has shrunk considerably in comparison to Figure~\ref{fig:introduction}. In fact, it looks like $h_\tau^+$ is very close to a time-shifted (and scaled) Gaussian when $\lvert \tau \rvert$ is large. 

It is therefore no surprise that time-shifting and scaling\footnote{For a precise explanation of how to time-shift and scale $h_\tau^\pm$ to obtain $f^\pm$, the reader is referred to Appendix~\ref{app:time-shift+scale}.} $h_\tau^\pm$, yields the functions
\begin{equation}
    \label{eq:counterexamples}
    f^\pm := f^\pm_\gamma := \varphi \pm \mathrm{i} \gamma \operatorname{T}_{1/a} \varphi, \qquad \gamma > 0,
\end{equation}
where the translation operator $\operatorname{T}_x$ is defined in equation~\eqref{eq:definitionoftranslation}. We emphasize that these functions satisfy 
\begin{equation*}
    \lim_{\gamma \to 0} f^\pm_\gamma = \varphi
\end{equation*}
pointwise and in $L^p(\R)$, for all $p \in [1,\infty]$. Therefore, there exist counterexamples to sampled Gabor phase retrieval which are arbitrarily close to the Gaussian.


Next, we rigorously show that $f^\pm$ are counterexamples to sampled Gabor phase retrieval by direct computation. We note that this already follows from \cite{alaifari2021phase,grohs2022foundational}. Since most of our proof will be used at a later point and for the convenience of the reader, we give an independent derivation here. In the following, we will use the well-known fact that the Gabor transform of the normalized Gaussian $\varphi$ is 
\begin{equation}\label{eq:gabortransformgaussian}
    \mathcal{G} \varphi (x,\omega) = 2^{1/4} \int_{\R} \varphi(t) \mathrm{e}^{-\pi(t-x)^2} \mathrm{e}^{-2\pi \mathrm{i} t \omega} \,\mathrm{d} t = \mathrm{e}^{-\pi\mathrm{i} x \omega} \mathrm{e}^{-\frac{\pi}{2} \left( x^2 + \omega^2 \right)},
 \end{equation}
 for $(x,\omega) \in \R^2$.

\begin{lemma}
    \label{lem:new_counterexamples}
    Let $a,\gamma > 0$ and let $f^\pm$ be defined as in equation~\eqref{eq:counterexamples}. Then, $f^+$ and $f^-$ do \emph{not} agree up to global phase and yet 
    \begin{equation*}
        \lvert \mathcal{G} f^+ \rvert = \lvert \mathcal{G} f^- \rvert \mbox{ on } \R \times a \mathbb{Z}.
    \end{equation*}
\end{lemma}

\begin{proof}
   The statement can be obtained by computing the Gabor transforms of $f^\pm$. By the linearity of the Gabor transform and the covariance property (cf.~\cite[Lemma 3.1.3 on p.~41]{groechenig2001foundations}), we find that 
    \begin{align*}
        \mathcal{G} f^\pm (x,\omega) &= \mathcal{G} \varphi (x,\omega) \pm \mathrm{i} \gamma \mathcal{G} \operatorname{T}_{1/a} \varphi (x,\omega)
        = \mathcal{G} \varphi (x,\omega) \pm \mathrm{i} \gamma \mathrm{e}^{-2\pi\mathrm{i} \frac{\omega}{a}} \mathcal{G} \varphi \left( x - \frac{1}{a},\omega \right) \\
        &= \mathrm{e}^{-\pi\mathrm{i} x \omega} \mathrm{e}^{-\frac{\pi}{2} \left( x^2 + \omega^2 \right)} \pm \mathrm{i} \gamma \mathrm{e}^{-2\pi\mathrm{i} \frac{\omega}{a}} \mathrm{e}^{-\pi\mathrm{i} \left(x - \frac{1}{a}\right) \omega} \mathrm{e}^{-\frac{\pi}{2} \left( \left(x - \frac{1}{a}\right)^2 + \omega^2 \right)} \\
        &= \mathrm{e}^{-\pi\mathrm{i} x \omega} \mathrm{e}^{-\frac{\pi}{2} \left( x^2 + \omega^2 \right)} \pm \mathrm{i} \gamma \mathrm{e}^{-\pi\mathrm{i} \left(x + \frac{1}{a}\right) \omega} \mathrm{e}^{-\frac{\pi}{2} \left( \left(x - \frac{1}{a}\right)^2 + \omega^2 \right)},
    \end{align*}
    for $(x,\omega) \in \R^2$. Therefore, we may compute
    \begin{align}
        \lvert \mathcal{G} f^\pm (x,\omega) \rvert &= \left\lvert \mathrm{e}^{-\frac{\pi}{2} \left( x^2 + \omega^2 \right)} \pm \mathrm{i} \gamma \mathrm{e}^{-\frac{\pi\mathrm{i} \omega}{a}} \mathrm{e}^{-\frac{\pi}{2} \left( \left(x - \frac{1}{a}\right)^2 + \omega^2 \right)} \right\rvert \nonumber \\
        &= \mathrm{e}^{-\frac{\pi}{2} \left( x^2 + \omega^2 \right)} \left\lvert 1 \pm \mathrm{i} \gamma \mathrm{e}^{\frac{\pi}{a}(x-\mathrm{i} \omega)} \mathrm{e}^{-\frac{\pi}{2a^2}} \right\rvert. \label{eq:Gaborfpm}
    \end{align}
    
    According to equation \eqref{eq:Gaborfpm}, the Gabor transform of $f^\pm$ is zero at $(x,\omega)$ if and only if 
    \begin{equation*}
        \mathrm{e}^{\frac{\pi}{a}(x-\mathrm{i} \omega)-\frac{\pi}{2a^2}} = \pm \frac{\mathrm{i}}{\gamma} = \mathrm{e}^{-\log \gamma \pm \frac{\pi \mathrm{i}}{2} + 2 \pi \mathrm{i} k}, 
    \end{equation*}
    for some $k \in \mathbb{Z}$, which is equivalent to 
    \begin{equation*}
        \frac{\pi}{a}(x-\mathrm{i} \omega) = \frac{\pi}{2a^2} - \log \gamma \pm \frac{\pi \mathrm{i}}{2} + 2 \pi \mathrm{i} k.
    \end{equation*}
    Therefore, the root sets of $\mathcal{G} f^\pm$ are given by 
    \begin{equation}
        \label{eq:root_set_Gcounterex}
        \left\{ \left( \frac{1}{2a} - \frac{a \log \gamma}{\pi}, \pm \frac{a}{2} + 2 a k \right) \,:\, k \in \mathbb{Z} \right\}.
    \end{equation}
    
    We note here that the root sets of $\mathcal{G} f^+$ and $\mathcal{G} f^-$ are different from each other so that $\mathcal{G} f^+$ and $\mathcal{G} f^-$ do not agree up to global phase. It follows by the linearity of the Gabor transform that $f^+$ and $f^-$ cannot agree up to global phase. Finally, we consider equation~\eqref{eq:Gaborfpm} once again to see that 
    \begin{align*}
        \lvert \mathcal{G} f^+ (x,a k) \rvert &= \mathrm{e}^{-\frac{\pi}{2} \left( x^2 + a^2 k^2 \right)} \left\lvert 1 + \mathrm{i} \gamma \mathrm{e}^{\frac{\pi}{a}(x- a \mathrm{i} k)} \mathrm{e}^{-\frac{\pi}{2a^2}} \right\rvert \\
        &= \mathrm{e}^{-\frac{\pi}{2} \left( x^2 + a^2 k^2 \right)} \left\lvert 1 + \mathrm{i} \gamma \mathrm{e}^{\frac{\pi x}{a}} \mathrm{e}^{- \pi \mathrm{i} k} \mathrm{e}^{-\frac{\pi}{2a^2}} \right\rvert \\
        &= \mathrm{e}^{-\frac{\pi}{2} \left( x^2 + a^2 k^2 \right)} \left\lvert 1 + \mathrm{i} (-1)^k \gamma \mathrm{e}^{\frac{\pi x}{a}} \mathrm{e}^{-\frac{\pi}{2a^2}} \right\rvert \\
        &= \mathrm{e}^{-\frac{\pi}{2} \left( x^2 + a^2 k^2 \right)} \left\lvert 1 - \mathrm{i} (-1)^k \gamma \mathrm{e}^{\frac{\pi x}{a}} \mathrm{e}^{-\frac{\pi}{2a^2}} \right\rvert \\
        &= \mathrm{e}^{-\frac{\pi}{2} \left( x^2 + a^2 k^2 \right)} \left\lvert 1 - \mathrm{i} \gamma \mathrm{e}^{\frac{\pi}{a}(x- a \mathrm{i} k)} \mathrm{e}^{-\frac{\pi}{2a^2}} \right\rvert \\
        &= \lvert \mathcal{G} f^- (x,a k) \rvert
    \end{align*}
    must hold, for $x \in \R$ and $k \in \mathbb{Z}$.
\end{proof}

\subsection{The set of counterexamples is dense in the space of square-integrable signals}\label{sec:classcounterexamplesdense}

Using the counterexamples $f^\pm$, which can be arbitrarily close to the Gaussian, we can show that if $\Lambda \subset \R^2$ is any set of equidistant parallel lines or any lattice, then \emph{the class of counterexamples $\mathfrak{C}(\Lambda)$ is dense in $L^2(\R)$}. Intuitively, our proof for this statement will work because the Bargmann transforms of the counterexamples $f^\pm$ are given by 
\begin{equation*}
    \mathcal{B} f^\pm (z) = 1 \pm \mathrm{i} \gamma \mathrm{e}^{-\tfrac{\pi}{2a^2} + \tfrac{\pi z}{a}}, \qquad z \in \bbC.
\end{equation*}
We can absorb the factor $\exp(-\pi/(2a^2))$ into $\gamma$ and obtain 
\begin{equation}\label{eq:Fock_multipliers}
    H^\pm_\delta(z) := 1 \pm \mathrm{i} \delta \mathrm{e}^{\tfrac{\pi z}{a}}, \qquad z \in \bbC,
\end{equation}
with $\delta>0$, which are entire functions of exponential type and thus in the Fock space. Additionally, $H^\pm_\delta$ converge to $1$ in $\mathcal{F}^2(\bbC)$ as $\delta \to 0$, do not agree up to global phase and satisfy $\lvert H^+_\delta \rvert = \lvert H^-_\delta \rvert$ on $\R + \mathrm{i} a \mathbb{Z}$: They are, in short, the ideal multipliers in the Fock space to transform a general function $f \in L^2(\R)$ into a ``close-by" counterexample.


\begin{theorem}\label{thm:classcounterexamplesdense}
    Let $a > 0$. Then, $\mathfrak{C}(\R \times a \Z)$ is dense in $L^2(\R)$.
\end{theorem}

\begin{proof}
    Let $\epsilon > 0$ and $f \in L^2(\R)$. We want to show that there exist $g^\pm \in L^2(\R)$ which do not agree up to global phase, are $\epsilon$-close to $f$ in $L^2(\R)$, i.e.~$\lVert f - g^\pm \rVert_2 < \epsilon$, and satisfy 
    \begin{equation*}
        \lvert \mathcal{G} g^+ \rvert = \lvert \mathcal{G} g^- \rvert \mbox{ on } \R \times a \Z.
    \end{equation*}

    According to \cite[Theorem 3.4.2 on p.~54]{groechenig2001foundations}, the monomials
    \begin{equation*}
        e_n(z) := \left(\frac{\pi^n}{n!}\right)^{1/2} z^n, \qquad n \in \mathbb{N}_0,~z \in \bbC,
    \end{equation*}
    form an orthonormal basis for the Fock space $\mathcal{F}^2(\bbC)$. Therefore, the space of complex polynomials is dense in the Fock space and we can find $P \in \bbC[z]$ such that 
    \begin{equation*}
        \lVert \mathcal{B} f - P \rVert_{\mathcal{F}} < \frac{\epsilon}{2}.
    \end{equation*}

    Let us now consider the functions $H^\pm_\delta$ defined in equation~\eqref{eq:Fock_multipliers} and note that $G_\delta^\pm := H_\delta^\pm \cdot P \in \mathcal{F}^2(\bbC)$ since $G_\delta^\pm$ are entire functions of exponential type. Hence, we can define $g^\pm_\delta := \mathcal{B}^{-1} G^\pm_\delta \in L^2(\R)$. To establish the necessary properties of $g^\pm_\delta$, we will work with their Bargmann transforms $G_\delta^\pm$. First, we note that 
    \begin{equation*}
        \lvert H_\delta^+(t + \mathrm{i} a k) \rvert = \lvert 1 + (-1)^k \mathrm{i} \delta \mathrm{e}^{\frac{\pi t}{a}} \rvert
        = \lvert 1 - (-1)^k \mathrm{i} \delta \mathrm{e}^{\frac{\pi t}{a}} \rvert = \lvert H_\delta^-(t + \mathrm{i} a k) \rvert,
    \end{equation*}
    for $t \in \R$ and $k \in \Z$. It follows that $\lvert G^+_\delta \rvert = \lvert G^-_\delta \rvert$ on $\R + \mathrm{i} a\Z$ and thus that $\lvert \mathcal{G} g^+_\delta \rvert = \lvert \mathcal{G} g^-_\delta \rvert$ on $\R \times a \Z$. Secondly, we note that $G^\pm_\delta$ do not agree up to global phase: Indeed, $H^\pm_\delta$ both have infinitely many roots but no root of $H^+_\delta$ is a root of $H^-_\delta$ and vice versa. At the same time, $P$ is a polynomial and has only finitely many roots. It follows that $G^+_\delta$ does have roots which are no roots of $G^-_\delta$ (and vice versa) and thus $G^+_\delta \not \sim G^-_\delta$. By the linearity of the Bargmann transform, we find $g^+_\delta \not\sim g^-_\delta$. Finally, we note that 
    \begin{equation*}
        \lVert P - P \cdot H^\pm_\delta \rVert_\mathcal{F} = \delta \lVert P \cdot \mathrm{e}^{\pi/a \cdot} \rVert_\mathcal{F}
    \end{equation*}
    and so there exists a $\delta > 0$ depending on $a$, $\epsilon$ and $P$ (which in turn depends on $f$ and $\epsilon$) such that 
    \begin{equation*}
        \lVert P - P \cdot H^\pm_\delta \rVert_\mathcal{F} < \frac{\epsilon}{2}.
    \end{equation*}
    We conclude that
    \begin{equation*}
        \lVert f - g^\pm_\delta \rVert_2 = \lVert \mathcal{B} f - H_\delta^\pm \cdot P \rVert_\mathcal{F} \leq \lVert \mathcal{B} f - P \rVert_\mathcal{F} + \lVert P - H_\delta^\pm \cdot P \rVert_\mathcal{F} < \epsilon.
    \end{equation*}
\end{proof}

\begin{remark}[Some explanations on the proof]
    \label{rem:interesting_examples}
    As $\mathcal{B} f \in \mathcal{F}^2(\bbC)$, for $f \in L^2(\R)$, we know that $\mathcal{B} f$ is either an entire function of exponential type or an entire function of second order. If $\mathcal{B} f$ is of second order, then its type is less or equal to $\pi/2$. 
   In the case that the type is strictly smaller than $\pi/2$, it holds that $\mathcal{B} f \cdot H^\pm_\delta \in \mathcal{F}^2(\bbC)$ and thus we can define 
    \begin{equation*}
        g^\pm_\delta := \mathcal{B}^{-1} \left( \mathcal{B} f \cdot H^\pm_\delta \right) \in L^2(\R),
    \end{equation*}
    with 
    \begin{equation*}
        \delta < \frac{\epsilon}{\lVert \mathcal{B} f \cdot \mathrm{e}^{\pi/a \cdot} \rVert_\mathcal{F}},
    \end{equation*}
    to obtain counterexamples which are $\epsilon$-close to $f$ in $L^2(\R)$.

   On the other hand, if $\mathcal{B} f$ is a second-order entire function of type $\pi/2$, it is not guaranteed that $\mathcal{B} f \cdot H^\pm_\delta$ is in the Fock space (see \cite{beneteau2010extremal} for two striking examples for why this can fail). As the only situation in which $\mathcal{B} f \cdot H^\pm_\delta$ may not be in the Fock space occurs when $\mathcal{B} f$ is exactly of order two and of type $\pi/2$, it seems obvious that the functions $f$ for which $\mathcal{B} f \cdot H^\pm_\delta \in \mathcal{F}^2(\bbC)$ holds must be dense in $L^2(\R)$. We can prove this by realizing that the complex polynomials are dense in $\mathcal{F}^2(\bbC)$. 

\end{remark}

We can adapt the argument above to work on a general set of parallel lines in the time-frequency plane. In order to do so, we use the functions 
\begin{equation*}
    H^\pm_\delta(z) := 1 \pm \mathrm{i} \delta \exp\left( \frac{\pi \mathrm{e}^{\mathrm{i} \theta}}{a} \left( z - \overline \lambda_0 \right) \right)
\end{equation*}
in our proof and note that the corresponding $g^\pm_\delta \in L^2(\R)$ satisfy 
\begin{equation*}
    \lvert \mathcal{G} g^+_\delta \rvert = \lvert \mathcal{G} g^-_\delta \rvert \mbox{ on } \operatorname{R}_\theta \left( \R \times a \Z \right) + \lambda_0,
\end{equation*}
where $a > 0$, $\lambda_0 \in \R^2 \simeq \bbC$, and $R_\theta : \R^2 \to \R^2$ denotes rotation by $\theta \in \R$ in $\R^2$ as defined in equation~\eqref{eq:definitionofrotation}. The statement for general lattices follows from the same consideration because all lattices are a subset of some set of infinitely many equidistant parallel lines. We therefore arrive at the following result.

\begin{theorem}
    \label{thm:density_counterexamples}
    Let $\Lambda \subset \R^2$ be a set of equidistant parallel lines or a lattice. Then, $\mathfrak{C}(\Lambda)$ is dense in $L^2(\R)$.
\end{theorem}

\begin{example}
    To illustrate our result, we create counterexamples to sampled Gabor phase retrieval which are close to the \emph{Hermite functions} $(H_n)_{n \geq 0} \in L^2(\R)$. To do so, we recall \cite[Theorem 3.4.2 on p.~54]{groechenig2001foundations} that
    \begin{equation*}
        \mathcal{B} H_n(z) = e_n(z) = \left(\frac{\pi^n}{n!}\right)^{1/2} z^n, \qquad z \in \bbC.
    \end{equation*}
    Therefore, the Gabor transform of the Hermite functions is 
    \begin{align*}
        \mathcal{G} H_n(x,\omega) &= \mathrm{e}^{-\pi\mathrm{i} x \omega} \mathcal{B} H_n(x - \mathrm{i} \omega) \mathrm{e}^{-\frac{\pi}{2}\left( x^2 + \omega^2 \right)} \\
        &= \left(\frac{\pi^n}{n!}\right)^{1/2} \mathrm{e}^{-\pi\mathrm{i} x \omega} \left( x - \mathrm{i} \omega \right)^n \mathrm{e}^{-\frac{\pi}{2}\left( x^2 + \omega^2 \right)},
    \end{align*}
    for $(x,\omega) \in \R^2$. Figure~\ref{fig:Gabor_Hermite} provides a plot of the Gabor transform (in magnitude) of the Hermite function for $n=5$.

    \begin{figure}
        \subfloat[$\lvert \mathcal{G} H_5 \rvert$]{\label{fig:Gabor_Hermite}
            \includegraphics[width=.49\textwidth]{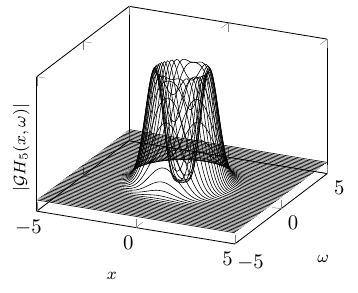}}
        \hfill
        \subfloat[$\lvert \mathcal{G} g_\delta^+ \rvert$, for $\delta = \tfrac{1}{50} \exp(-10 \pi)$.]{\label{fig:counterexample}
            \includegraphics[width=.49\textwidth]{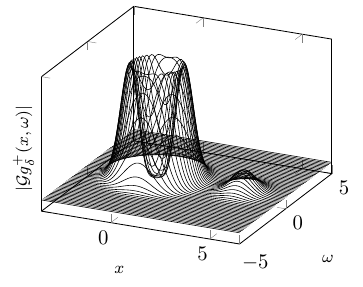}}
        \caption{The Gabor transform magnitude of the fifth Hermite function (Fig.~\ref{fig:Gabor_Hermite}) and of a counterexample $g_\delta^+$ to sampled Gabor phase retrieval on $\R \times \tfrac{1}{4} \Z$ (Fig.~\ref{fig:counterexample}).}
        \label{fig:somename}
    \end{figure}

    To find a counterexample which is close to $H_n$, we can define 
    $$g^\pm_\delta := \mathcal{B}^{-1} ( \mathcal{B} H_n \cdot H^\pm_\delta )$$ 
    as described in Remark~\ref{rem:interesting_examples}. For visualization purposes, we are interested in the spectrogram of $g^\pm_\delta$ and compute
    \begin{align*}
        \mathcal{G} g^\pm_\delta(x,\omega) &= \mathrm{e}^{-\pi\mathrm{i} x \omega} \mathcal{B} g^\pm_\delta(x - \mathrm{i} \omega) \mathrm{e}^{-\frac{\pi}{2}\left( x^2 + \omega^2 \right)} \\
        &= \mathrm{e}^{-\pi\mathrm{i} x \omega} \mathcal{B} H_n(x - \mathrm{i} \omega) \mathrm{e}^{-\frac{\pi}{2}\left( x^2 + \omega^2 \right)} \cdot H^\pm_\delta(x - \mathrm{i} \omega) \\
        &= \mathcal{G} H_n(x,\omega) \cdot H^\pm_\delta(x - \mathrm{i} \omega) \\
        &= \left(\frac{\pi^n}{n!}\right)^{1/2} \mathrm{e}^{-\pi\mathrm{i} x \omega} \left( x - \mathrm{i} \omega \right)^n \left( 1 \pm \mathrm{i} \delta \mathrm{e}^{\frac{\pi (x- \mathrm{i} \omega)}{a}} \right) \mathrm{e}^{-\frac{\pi}{2}\left( x^2 + \omega^2 \right)}.
    \end{align*}
    For comparison, we plot the magnitude of the above Gabor transform (for $n = 5$, $a = \tfrac{1}{4}$ and $\delta = \tfrac{1}{50} \exp(-10 \pi)$) in Figure \ref{fig:counterexample}.
\end{example}

\begin{remark}
    It is notable that if $f \in L^2(\R)$ is such that $\mathcal{B} f$ is \emph{not} of second order and of type $\pi/2$, then there are counterexamples of the form
    \begin{equation*}
        g^\pm_\delta = f \pm \mathrm{i} \delta \mathcal{B}^{-1} \left( z \mapsto \mathrm{e}^{\frac{\pi z}{a}} \mathcal{B} f (z) \right);
    \end{equation*}
    i.e.~$g^\pm_\delta$ are small additive perturbations of our original signals $f$.
\end{remark}

\subsection{Recovering the Gaussian from Gabor transform magnitudes  on a lattice}\label{ssec:uniqueness_for_Gaussian}

Theorem~\ref{thm:density_counterexamples} states that the counterexamples $\mathfrak{C}(\Lambda)$ are dense in $L^2(\R)$; in other words, for every $f \in L^2(\R)$, there is a counterexample to sampled Gabor phase retrieval which is arbitrarily close to $f$. Naturally, one might wonder whether all $f \in L^2(\R) \setminus \{0\}$ are counterexamples. In this subsection, we show that this is not true and that, in particular, \emph{the normalized Gaussian can be recovered from Gabor magnitude samples on sufficiently fine square lattices}. 

In order to accomplish this, we relate sampled Gabor phase retrieval to a sampling problem in the Fock space of entire functions. In this way, our result reduces to identifying all $F \in \mathcal{F}^2(\mathbb{C})$ which have unit absolute value on a square lattice. This reduction relates our problem to some form of maximum modulus (or Phragm{\'e}n--Lindel{\"o}f) principle (cf.~\cite[Section~5.1 on pp.~165--168 and Section~5.6 on pp.~176--181]{titchmarsh1939theory}): Indeed, we are considering a second order entire function $F$ which is bounded on all lattice points; this suggests that $F$ should be constant in the entire complex plane as long as the lattice is dense enough. This intuition is correct and follows from an elegant result discovered independently by V.~Ganapathy Iyer \cite{ganapathy1936note} and Albert Pfluger \cite{pfluger1937analytic} in 1936.

\begin{theorem}[{Cf.~\cite[Theorem I A on p.~305]{pfluger1937analytic}}]
    \label{thm:coreresult}
    Let $h$ be an entire function such that 
    \[
        \limsup_{r \to \infty} \frac{\log M_h(r)}{r^2} < \frac{\pi}{2},
    \]
    where $M_h(r) := \max_{\lvert z \rvert = r} \lvert h(z) \rvert$. If there exists 
    a constant $\kappa > 0$ such that 
    \[
        \lvert h(m+\mathrm{i} n) \rvert \leq \kappa, \qquad m,n \in \mathbb{Z},
    \]
    then $h$ is constant.
\end{theorem}

We can now prove the following result.

\begin{theorem}
    \label{thm:original}
    Let $0 < a < 1$ and $f \in L^2(\R)$ be such that 
    \[
        \lvert \mathcal{G} f(x,\omega) \rvert^2  = \mathrm{e}^{-\pi \left( x^2 + \omega^2 \right)}  = \lvert \mathcal{G} \varphi(x,\omega) \rvert^2, \qquad (x,\omega) \in a \mathbb{Z}^2.
    \]
    Then, there exists an $\alpha \in \R$ such that $f = \mathrm{e}^{\mathrm{i} \alpha} \varphi$.
\end{theorem}

\begin{proof}
    Let us consider the entire function $h(z) := \mathcal{B}f(az)$ for $z \in \mathbb{C}$. We directly estimate
    \[
        \lvert h(z) \rvert = \lvert \mathcal{B}f(az) \rvert \leq \lVert \mathcal{B}f \rVert_\mathcal{F} \cdot  \mathrm{e}^{\frac{\pi}{2} \lvert az \rvert^2} 
        = \lVert f \rVert_2 \cdot \mathrm{e}^{\frac{\pi a^2}{2} \lvert z \rvert^2}, \qquad z \in \mathbb{C},
    \]
    following \cite[Proposition~3.4.1 and Theorem~3.4.2 on p.~54]{groechenig2001foundations} which shows that
    \[
        \limsup_{r \to \infty} \frac{\log M_h(r)}{r^2} \leq \limsup_{r \to \infty} \left(\frac{\log \lVert f \rVert_2}{r^2} + \frac{\pi a^2}{2} \right) = \frac{\pi a^2}{2} < \frac{\pi}{2}.
    \]

    Additionally, equation~\eqref{eq:Gabor_Bargmann} and our assumption on $f$ imply that 
    \begin{equation*}
        \lvert h(m + \mathrm{i} n) \rvert = \lvert \mathcal{B} f (am + \mathrm{i} a n) \rvert = \lvert \mathcal{G} f (am,-an) \rvert \mathrm{e}^{\frac{\pi a^2}{2}(m^2 + n^2)} = 1.
    \end{equation*}
    Therefore, Theorem~\ref{thm:coreresult} shows that $h$ is constant which together with $\lvert h \rvert = 1$ implies that there exists an $\alpha \in \R$ such that $h = \mathrm{e}^{\mathrm{i} \alpha}$. Since the Bargmann transform is unitary \cite[Theorem~3.4.3 on p.~56]{groechenig2001foundations} and linear, and since $\mathcal{B} \varphi = 1$, it follows that $f = \mathrm{e}^{\mathrm{i} \alpha} \varphi$ as desired.
\end{proof}

\begin{remark}
    A natural confusion that might arise in connection with Theorem~\ref{thm:original} is in how far it is different from the result in \cite{grohs2023injectivity} on shift-invariant spaces with Gaussian generator,
    \[
        V_\beta^p(\varphi) := \left\{ f \in L^p(\R) \,:\, f = \sum_{k \in \mathbb{Z}} c_k \operatorname{T}_{\beta k} \varphi,~c \in \ell^p(\mathbb{Z}) \right\}, 
    \]
    where $p \in [1,\infty]$ and $\beta \in (0,\infty)$. The aforementioned result applied to the Gaussian states that if $\beta > 0$ and $0 < a < \beta/2$ are such that $a \beta \not\in \mathbb{Q}$, then the only functions $f \in V_\beta^1(\varphi)$ satisfying 
    \begin{equation*}
        \lvert \mathcal{G} f(x,\omega) \rvert^2 = \lvert \mathcal{G} \varphi(x,\omega) \rvert^2 = \mathrm{e}^{-\pi \left( x^2 + \omega^2 \right)}, \qquad (x,\omega) \in a \mathbb{Z}^2,
    \end{equation*}
    are of the form $f = \mathrm{e}^{\mathrm{i} \alpha} \varphi$, where $\alpha \in \R$. The difference is the assumption $f \in V_\beta^1(\varphi)$ which stands in contrast to the weaker assumption $f \in L^2(\R)$ in Theorem~\ref{thm:original}. In short, Theorem~\ref{thm:original} implies that the Gaussian can be distinguished from all other functions in $L^2(\R)$ by looking at its sampled Gabor transform magnitude measurements while \cite[Theorem 3.10 on p.~188]{grohs2023injectivity} only implies that it is distinguishable from the functions in $V_\beta^1(\varphi) \subsetneq L^2(\R)$.
\end{remark}

\section{On the stability of Gabor phase retrieval}\label{sec:stability}

Having discussed uniqueness from samples in the previous section, we now turn to the question of stability and the lack of connection between these two properties. Let $\Omega\subseteq\R^2$ be a domain and let $1\leq p\leq\infty$.  
 Given $f\in M^{p}(\R)$, we denote by $|\mathcal{G}f_{|_{\Omega}}|$ the magnitude of the Gabor transform of $f$ on $\Omega$. We are interested in the local Lipschitz constant of $f$ on $\Omega$, i.e. the smallest constant $C>0$ such that
\begin{equation}\label{eq:local_stability_constant}
    \inf_{\alpha\in\R}\|\mathcal{G}f-e^{i\alpha}\mathcal{G}g\|_{L^p(\Omega)} 
    \leq C \||\mathcal{G}f_{|_{\Omega}}|-|\mathcal{G}g_{|_{\Omega}}|\|_{\mathfrak{B}}, \qquad \mbox{for all } g \in M^{p}(\R).
\end{equation}
Here, $\|\cdot\|_{\mathfrak{B}}$ denotes the norm of a Banach space $\mathfrak{B}$ in which the space of measurements $|\mathcal{G}f_{|_{\Omega}}|$, for $f\in M^{p}(\R)$, lie. In particular, we denote the local Lipschitz constant by $c_{p,\Omega}(f)$ when $\|\cdot\|_{\mathfrak{B}}$ is given by \eqref{eq:normmeasuremnts}, i.e. $c_{p,\Omega}(f)$ is the best possible constant $C>0$ for which 
\begin{equation}\label{eq:local_stability_constant2}
    \inf_{\alpha\in\R}\|\mathcal{G}f-e^{i\alpha}\mathcal{G}g\|_{L^p(\Omega)} 
    \leq C \||\mathcal{G}f_{|_{\Omega}}|-|\mathcal{G}g_{|_{\Omega}}|\|_{\mathcal{D}_{p,q}^{1,4}(\Omega)}, \quad \mbox{for all } g \in M^{p}(\R).
\end{equation}

A large constant $c_{p,\Omega}(f)$ indicates that the problem of recovering $\mathcal{G} f_{|_{\Omega}}$ from $|\mathcal{G}f_{|_{\Omega}}|$ cannot be controlled well since there exists a function $g\in M^{p}(\R)$ with $|\mathcal{G}g_{|_{\Omega}}|$ very close to $|\mathcal{G}f_{|_{\Omega}}|$ while the distance between $\mathcal{G}f_{|_{\Omega}}$ and $\mathcal{G}g_{|_{\Omega}}$ is not small. Consequently, the problem of recovering $f$ from $|\mathcal{G}f_{|_{\Omega}}|$ is not well controlled either since  
\[
    \inf_{\alpha\in\R}\|f-e^{i\alpha}g\|_{M^{p}(\R)}
    = \inf_{\alpha\in\R}\|\mathcal{G}f-e^{i\alpha}\mathcal{G}g\|_{L^p(\R^2)}
    \geq \inf_{\alpha\in\R}\|\mathcal{G}f-e^{i\alpha}\mathcal{G}g\|_{L^p(\Omega)}.
\]
On the other hand, a small $c_{p,\Omega}(f)$ translates into good stability guarantees for the recovery of $\mathcal{G} f_{|_{\Omega}}$ from $|\mathcal{G}f_{|_{\Omega}}|$. We observe that, if $\Omega\subsetneq\R^2$, this does not guarantee that the problem of recovering $f$ from $|\mathcal{G}f_{|_{\Omega}}|$ is stable. However, if we suppose that $f$ is $\epsilon$-concentrated on $\Omega$, i.e.~$f$ satisfies
\begin{equation}\label{eq:epsilonconcentrated}
\|\mathcal{G}f\|_{L^p(\R^2\setminus\Omega)}\leq\epsilon,    
\end{equation}
for some small $\epsilon>0$, then we obtain a weaker notion of stability for the recovery of $f$ from $|\mathcal{G}f_{|_{\Omega}}|$ in the sense that 
$$
    \inf_{\alpha\in\R}\|f-e^{i\alpha}g\|_{M^{p}(\R)}=
     \inf_{\alpha\in\R}\|\mathcal{G}f-e^{i\alpha}\mathcal{G}g\|_{L^p(\R^2)}\leq c_{p,\Omega}(f) \||\mathcal{G}f_{|_{\Omega}}|-|\mathcal{G}g_{|_{\Omega}}|\|_{\mathfrak{B}}+2\epsilon,
$$
for any $g\in M^p(\R)$ that is $\epsilon$-concentrated on $\Omega$. 



As discussed in Section~\ref{ssec:prelim_LaplacePoincareCheeger}, the stability constant $c_{p,\Omega}(f)$ can be controlled by the weighted Poincaré constant: 
the smaller the Poincar\'e constant, the better the local stability of Gabor phase retrieval at $f$ (cf.~equation~(\ref{eq:stab-GR})). The relation to the Cheeger constant is inversely proportional. More precisely, the overall picture is: 
\[
    c_{p,\Omega}(f) \lesssim C_{\mathrm{poinc}}(p,\Omega,|\mathcal{G}f|^p) \lesssim h_{p,\Omega}(f)^{-1}.
\]

In what follows, we will show that the construction of the counterexamples \eqref{eq:counterexamples} leads to the existence of function perturbations which break uniqueness from samples while not affecting stability. Hence, much to our surprise, the answer to Question~\ref{conjecture} is negative. Note that the precise statement in our Theorem \ref{thm:main} slightly differs from Question~\ref{conjecture} and we discuss these technicalities in Section~\ref{ssec:falseconjecture}. 

\subsection{On the variation of the weighted Poincar\'e constant}\label{ssec:variation}

A natural question concerning weighted Poincar\'e inequalities is how the Poincar\'e constant $C_{\mathrm{poinc}}(p,\Omega,w)$ changes under variations of the weight $w$. Lemma~\ref{lem:var_PWconst} provides a simple result in that direction. We first state a classical fact that we exploit in its proof. 

\begin{lemma}[\cite{dyda2013weighted}]
    \label{lem:usefullem}
    Let $1 \leq p < \infty$, let $\Omega \subset \mathbb{R}^2$ be a domain, and let $w$ be a weight on $\Omega$. Then, for every $F \in L^p(\Omega,w)$, it holds that 
    \begin{equation}\label{eq:usefuleq}
      \inf_{c \in \mathbb{R}} \lVert F - c \rVert_{L^p(\Omega,w)} \leq \|F-F_{\Omega}^{w}\|_{L^p(\Omega,w)} \leq 2 \inf_{c \in \mathbb{R}} \lVert F - c \rVert_{L^p(\Omega,w)}.
    \end{equation}
\end{lemma}

\begin{lemma}
    \label{lem:var_PWconst}
    Let $1 \leq p < \infty$, let $\Omega \subset \mathbb{R}^2$ be a domain, and let $w$ be a weight on $\Omega$. 
    Let $w'$ be another weight on $\Omega$ which satisfies
    \begin{equation}\label{eq:ratioweights}
        A w(x) \leq w'(x) \leq B w(x), \qquad x \in \Omega,
    \end{equation}
    for some constants $0 < A \leq B < \infty$. Then, 
    it holds that 
    \begin{equation}\label{eq:ratiopoincareconstants}
        \frac{A^{1/p}}{ 2B^{1/p}}\  C_{\mathrm{poinc}}(p,\Omega,w) \leq C_{\mathrm{poinc}}(p,\Omega,w') \leq \frac{2 B^{1/p}}{A^{1/p}}\  C_{\mathrm{poinc}}(p,\Omega,w).
    \end{equation}
\end{lemma}

\begin{proof}
By Equation~\eqref{eq:ratioweights}, it follows that for every $F \in L^p(\Omega,w)$, the integral inequality
    \begin{equation*}
        A \int_\Omega \lvert F(x) \rvert^p w(x) \,\mathrm{d} x \leq  \int_\Omega \lvert F(x) \rvert^p w'(x) \,\mathrm{d} x \leq B \int_\Omega \lvert F(x) \rvert^p w(x) \,\mathrm{d} x
    \end{equation*}
holds true. Hence, we have that 
    \begin{equation*}
        A^{1/p} \lVert F \rVert_{L^p(\Omega,w)} \leq \lVert F \rVert_{L^p(\Omega,w')} \leq B^{1/p} \lVert F \rVert_{L^p(\Omega,w)}
    \end{equation*}
    and consequently 
    $L^p(\Omega,w) = L^p(\Omega,w')$ as well as $W^{1,p}(\Omega,w) = W^{1,p}(\Omega,w')$. In particular, we obtain that for every $F \in L^p(\Omega,w)$ and for every $c\in\mathbb{R}$,
        \begin{equation*}
        A^{1/p} \lVert F-c \rVert_{L^p(\Omega,w)} \leq \lVert F-c \rVert_{L^p(\Omega,w')} \leq B^{1/p} \lVert F-c \rVert_{L^p(\Omega,w)},
    \end{equation*}
    as well as
    \begin{equation*}
        A^{1/p} \lVert \nabla F \rVert_{L^p(\Omega,w)} \leq \lVert \nabla F \rVert_{L^p(\Omega,w')} \leq B^{1/p} \lVert \nabla F \rVert_{L^p(\Omega,w)}.
    \end{equation*}
    
    We can now prove the upper and lower bounds in equation~\eqref{eq:ratiopoincareconstants}. By using Lemma~\ref{lem:usefullem} and equation~\eqref{eq:poincdefnsup} along with the above inequalities, we find that 
    \begin{align*}
        C_{\mathrm{poinc}}(p,\Omega,w') &\begin{multlined}[t]
            \leq 2 \sup \bigg\{ \inf_{c \in \mathbb{R}} \frac{\lVert F - c \rVert_{L^p(\Omega,w')}}{\lVert \nabla F \rVert_{L^p(\Omega,w')}} \,: \\
            \qquad F \in W^{1,p}(\Omega,w') \cap \mathcal{M}(\Omega),~F \neq \mathrm{const.} \bigg\}
        \end{multlined} \\
        &\begin{multlined}[t]
            \leq \frac{2 B^{1/p}}{A^{1/p}} \sup \bigg\{ \inf_{c \in \mathbb{R}} \frac{\lVert F - c \rVert_{L^p(\Omega,w)}}{\lVert \nabla F \rVert_{L^p(\Omega,w)}} \,:\\
            \qquad F \in W^{1,p}(\Omega,w)\cap \mathcal{M}(\Omega),~F \neq \mathrm{const.} \bigg\} 
        \end{multlined} \\
        &\begin{multlined}[t]
            \leq \frac{2 B^{1/p}}{A^{1/p}} \sup \bigg\{ \frac{\lVert F - F_\Omega^{w} \rVert_{L^p(\Omega,w)}}{\lVert \nabla F \rVert_{L^p(\Omega,w)}} \,: \\
            \qquad F \in W^{1,p}(\Omega,w)\cap \mathcal{M}(\Omega),~F \neq \mathrm{const.} \bigg\}
        \end{multlined} \\
        &= \frac{2 B^{1/p}}{A^{1/p}} \ C_{\mathrm{poinc}}(p,\Omega,w).
    \end{align*}
    
    We can essentially repeat the same argument to show the lower bound
    \begin{align*}
        C_{\mathrm{poinc}}(p,\Omega,w') &\geq \frac{A^{1/p}}{2 B^{1/p}} \ C_{\mathrm{poinc}}(p,\Omega,w)
    \end{align*}
    which concludes the proof.
\end{proof}

\subsection{Answering  Question~\ref{conjecture}}\label{ssec:falseconjecture}

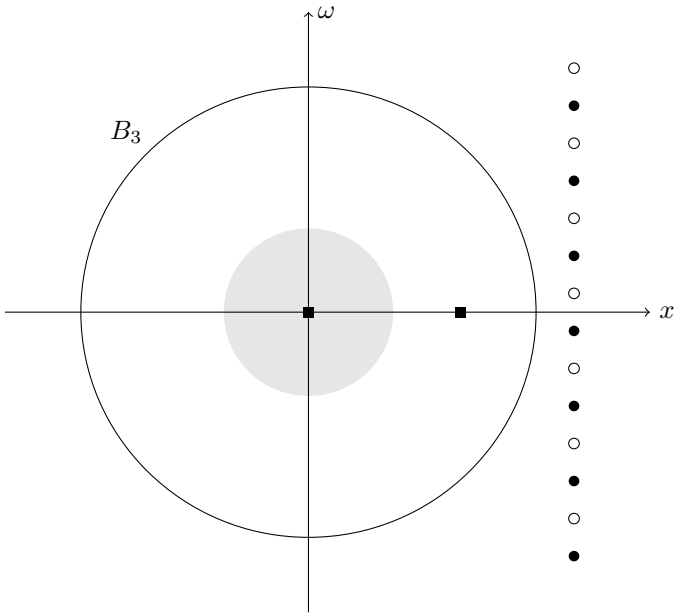
\begin{figure}
    \centering
    \begin{tikzpicture}

        \fill[gray!20] (0,0) circle (1.1166);
        \draw (0,0) circle (3);
        \node at (-2.4,2.4) {$B_3$};

        \draw[->] (-4,0)--(4.5,0) node[right] {$x$};
        \draw[->] (0,-4)--(0,4) node[right] {$\omega$};

        \fill[black] (0,0) +(-2pt,-2pt) rectangle +(2pt,2pt);
        \fill[black] (2,0) +(-2pt,-2pt) rectangle +(2pt,2pt);

        \draw (3.5,0.25) circle (2pt);
        \fill[black] (3.5,-0.25) circle (2pt);
        \foreach \k in {1,...,3}
        {
        \draw (3.5,{1/4+\k}) circle (2pt);
        \fill[black] (3.5,{-0.25+\k}) circle (2pt);
        \draw (3.5,{1/4-\k}) circle (2pt);
        \fill[black] (3.5,{-0.25-\k}) circle (2pt);
        }
        
        \draw [black] plot [only marks, mark size=2, mark=square*] coordinates {(0,0) (2,0)};
    \end{tikzpicture}
    \caption{We consider $a = 1/2$ and $\gamma = \exp(-5\pi)$. The roots of $\mathcal{G} f^+$ are indicated by circles and the roots of $\mathcal{G} f^-$ are indicated by disks. We have also drawn the local maxima of $\mathcal{G} f^\pm$ as squares and indicated the region on which 99\% of the $L^2$-mass of $\mathcal{G} f^\pm$ is concentrated in light gray. We highlight that we have chosen $\gamma < \gamma_0(1/2,R=3) = \exp(-4\pi)$ such that the roots of $\mathcal{G} f^\pm$ fall outside the open ball of radius $R=3$. We also note that there is no gray region around the local maximum at $(2,0)$ indicating that very little mass is concentrated on the small bump.}
    \label{fig:intersecting_splitting_roots}
\end{figure}

We apply Lemma~\ref{lem:var_PWconst} to the special case where the weights are given by
$$
w=|\mathcal{G}\varphi|^p,\quad  w_\pm'=|\mathcal{G}f_{\pm}|^p,
$$
for some $1 \leq p < \infty$, and where $f_{\pm}$ denotes the counterexamples  
\begin{equation*}
    f_\pm = \varphi \pm \mathrm{i} \gamma \operatorname{T}_{1/a} \varphi, \qquad \gamma > 0,\, a>0,
\end{equation*}
constructed in Section~\ref{ssec:counterexamples}. By equations~\eqref{eq:gabortransformgaussian} and \eqref{eq:Gaborfpm}, we have that for all $(x,\omega) \in \R^2$,
\begin{equation*}
\lvert \mathcal{G} f_\pm (x,\omega) \rvert^p = \lvert \mathcal{G} \varphi (x,\omega) \rvert^p \left\lvert 1 \pm \mathrm{i} \gamma \mathrm{e}^{\frac{\pi}{a}(x-\mathrm{i} \omega)} \mathrm{e}^{-\frac{\pi}{2a^2}} \right\rvert^p.
\end{equation*}
    The expression of the Gabor transform of the counterexamples $f_\pm$ allows us to show that for all $a,R > 0$, there exists a $\gamma_0 = \gamma_0(a,R) > 0$ such that for all $\gamma \in (0,\gamma_0)$, the roots of $\mathcal{G} f_\pm$ fall outside of the cube $[-R,R]^2$. By equation~\eqref{eq:gabortransformgaussian}, $\lvert \mathcal{G} \varphi (x,\omega) \rvert>0$ for all $(x,\omega)\in\R^2$ and, thus, the roots of $\mathcal{G} f_\pm$ correspond to the roots of $(x,\omega) \mapsto \lvert 1 \pm \mathrm{i} \gamma \mathrm{e}^{\frac{\pi}{a}(x-\mathrm{i} \omega)} \mathrm{e}^{-\frac{\pi}{2a^2}} \rvert$. Furthermore, by the reverse triangular inequality 
    \[
        \left\lvert 1 \pm \mathrm{i} \gamma \mathrm{e}^{\frac{\pi}{a}(x-\mathrm{i} \omega)} \mathrm{e}^{-\frac{\pi}{2a^2}} \right\rvert\geq 1- \gamma \left\lvert\mathrm{e}^{\frac{\pi}{a}(x-\mathrm{i} \omega)} \mathrm{e}^{-\frac{\pi}{2a^2}} \right\rvert,
    \]
    and the continuous function $(x,\omega) \mapsto \lvert \mathrm{e}^{\frac{\pi}{a}(x-\mathrm{i} \omega)} \mathrm{e}^{-\frac{\pi}{2a^2}} \rvert$ attains its maximum on $[-R,R]^2$. Hence, the Gabor transforms $\mathcal{G} f_\pm$ have no roots in $[-R,R]^2$ provided that we choose
    \[
        \gamma < \left( \max_{x,\omega \in [-R,R]} \left\lvert \mathrm{e}^{\frac{\pi}{a}(x-\mathrm{i} \omega)} \mathrm{e}^{-\frac{\pi}{2a^2}} \right\rvert \right)^{-1}.
    \]
 We can also determine $\gamma_0$ precisely by examining the root sets of $\mathcal{G} f^\pm$ which are 
    \begin{equation*}
        \left\{ \left( \frac{1}{2a} - \frac{a \log \gamma}{\pi}, \pm \frac{a}{2} + 2 a k \right) \,:\, k \in \mathbb{Z} \right\},
    \end{equation*}
    as proven in equation \eqref{eq:root_set_Gcounterex}. (We have visualized these roots in Figure~\ref{fig:intersecting_splitting_roots}.) Indeed, let us consider arbitrary but fixed $a,R > 0$ and set 
    \begin{equation*}
        \gamma_0 := \mathrm{e}^{-\frac{\pi}{a}\left( R - \frac{1}{2a} \right)}.
    \end{equation*}
    Then, it holds that for all $\gamma \in (0,\gamma_0]$, the roots of $\mathcal{G} f_\pm$ fall outside the strip $(-R,R) \times \R$.
Hence, we restrict to a bounded domain $\Omega\subseteq\R^2$ and choose $R>0$ such that $\Omega\subseteq (-R,R)\times\R$. If $\gamma < \gamma_0$ as above, all the roots of $|\mathcal{G}f_{\pm}|^p$ fall outside of $\Omega$. As a consequence of the extreme value theorem, there exist $0<A_{\gamma}\leq B_{\gamma}<\infty$ such that 
\[
A_{\gamma}\leq \left\lvert 1 \pm \mathrm{i} \gamma \mathrm{e}^{\frac{\pi}{a}(x-\mathrm{i} \omega)} \mathrm{e}^{-\frac{\pi}{2a^2}} \right\rvert^p \leq B_{\gamma},\qquad (x,\omega) \in \Omega.
\]
More precisely, given any $0<\delta<1$, the stronger condition $$\gamma<\delta\mathrm{e}^{-\frac{\pi}{a}\left( R - \frac{1}{2a} \right)}$$ implies that for all $(x,\omega) \in \Omega$,
\[
 (1-\delta)^{p} \leq \left\lvert 1 \pm \mathrm{i} \gamma \mathrm{e}^{\frac{\pi}{a}(x-\mathrm{i} \omega)} \mathrm{e}^{-\frac{\pi}{2a^2}} \right\rvert^p \leq  (1+\delta)^p,
\]
and consequently 
\begin{equation}\label{eq:boundingweights}
   (1-\delta)^{p} \ \lvert \mathcal{G} \varphi (x,\omega) \rvert^p \leq \lvert \mathcal{G} f_\pm (x,\omega) \rvert^p \leq (1+\delta)^p \lvert \mathcal{G} \varphi (x,\omega) \rvert^p,
\end{equation}
for all $(x,\omega) \in \Omega$.

\begin{corollary}
\label{cor:maincor}
    Let $1 \leq p < \infty$, let $\Omega \subset \mathbb{R}^2$ be a bounded domain and let $a > 0$. Then, for any $0<\delta<1$, there exists a constant $\gamma_\delta = \gamma_\delta (a,\Omega) > 0$ such that for all $\gamma<\gamma_\delta$, it holds that 
    \begin{align*}
        \frac{(1-\delta)}{2(1+\delta)} \cdot C_{\mathrm{poinc}}\left(p,\Omega,\lvert \mathcal{G} \varphi \rvert^p \right) &\leq C_{\mathrm{poinc}}\left(p,\Omega,\lvert \mathcal{G} f_\pm \rvert^p\right) \\
        &\leq \frac{2(1+\delta)}{(1-\delta)} \cdot C_{\mathrm{poinc}}\left(p,\Omega,\lvert \mathcal{G} \varphi \rvert^p \right).
    \end{align*}
\end{corollary}

\begin{proof}
    The proof follows by applying Lemma~\ref{lem:var_PWconst} along with equation~\eqref{eq:boundingweights}.
\end{proof}

\begin{remark}
Theorem~B.7 together with Theorem~B.8 in \cite{grohs2019stable} ensure that 
$$C_{\mathrm{poinc}}\left(p,\Omega,\lvert \mathcal{G} \varphi \rvert^p \right)<\infty$$
whenever $p\in[1,2]$ and $\Omega\subseteq\R^2$ is a bounded domain with Lipschitz boundary.
\end{remark}

Let $\nu>0$ and let $B_R$ denote the ball of radius $R > 0$ centered at $0$. We recall the notation $\mathcal{M}_{\nu}(B_R)$ introduced in Question~\ref{conjecture} for the class of functions
\[
\mathcal{M}_{\nu}(B_R)=\{f\in M^{p}(\R): c_{p,B_R}(f)\leq \nu\}.
\]
The following theorem is our main result. It provides a theoretical foundation for our claim that the signal class $\mathcal{M}_\nu(\R^2)$ cannot serve as a prior for uniqueness in sampled Gabor phase retrieval. More precisely, it states that there exists $\nu > 0$ for which the signal class $\mathcal{M}_\nu(B_R)$, $R > 1$, contains functions that do not agree up to global phase but whose Gabor transform magnitudes agree on a rectangular lattice $\Lambda$ --- no matter how large we choose $R > 1$ and how fine we choose the lattice. This clearly implies the existence of counterexamples for sampled Gabor phase
retrieval on arbitrary rectangular lattices in all the classes of functions $\mathcal{M}_{\overline{\nu}}(B_R)$, $\overline{\nu}>\nu$, $R>1$. Observe that every rectangular lattice $\Lambda$ is contained in a set of parallel lines; that is, there exists $a>0$ such that $\Lambda\subseteq\R\times a\Z$.

\begin{theorem}[Main result]\label{thm:main}
Let $p\in[1,2)$, $q\in(2p/(2-p),\infty)$. There exists $\nu>0$ such that, for all $R>1$ and for all $a>0$, there exist $f,g\in \mathcal{M}_{\nu}(B_R)$ such that $f \not\sim g$ but
\[
    |\mathcal{G}f(x,\omega)|=|\mathcal{G}g(x,\omega)|,\quad (x,\omega)\in \R\times a\Z.
\]
\end{theorem}

\begin{proof}
We show that there exists $\nu>0$ such that, for all $R>1$ and for all $a>0$, there exists $\gamma>0$ such that 
\[
c_{p,B_R}(f_{\pm})\leq \nu,
\]
where 
\[
f_\pm = \varphi\pm\mathrm{i}\gamma \operatorname{T}_{\frac{1}{a}}\varphi.
\]
We already know from Section~\ref{ssec:counterexamples} that $f_\pm \in M^{p}(\R)$, $f_+\not\sim f_-$ and 
\[
    |\mathcal{G}f_+(x,\omega)|=|\mathcal{G}f_-(x,\omega)|,\quad (x,\omega)\in \R\times a\mathbb{Z}.
\]
Let $R>1$ and $a>0$. We choose 
\[
    \gamma < \delta \cdot \min\left\{1, \mathrm{e}^{-\frac{\pi}{a}\left( R - \frac{1}{2a} \right)} \right\},
\]
with $0<\delta<1$. The condition 
$\gamma<\delta\mathrm{e}^{-\frac{\pi}{a}\left( R - \frac{1}{2a} \right)}$ ensures that the roots of $f_\pm$ fall outside the ball $B_R$. Theorem~5.9 in \cite{grohs2019stable} states that 
\begin{equation}\label{eq:boundingstabilityconstant}
c_{p,B_R}(f_{\pm})\leq c(1+C_{\mathrm{poinc}}(p,B_R,|\mathcal{G}f_{\pm}|^p)),
\end{equation}
where $c>0$ is a constant depending on $p$, $q$ and monotonically increasingly on
\begin{equation}\label{eq:max}
 \max\{\|\mathcal{G}f_{\pm}\|_{L^p(B_R)}/\|\mathcal{G}f_{\pm}\|_{L^\infty(B_R)},\|\mathcal{V}_{\varphi'}f_{\pm}\|_{L^\infty(B_R)}/\|\mathcal{G}f_{\pm}\|_{L^\infty(B_R)}\},  
\end{equation}
where $\varphi'$ denotes the first derivative of the Gaussian $\varphi$.

By Corollary~\ref{cor:maincor}, we know how to upper bound the weighted Poincar\'e constant in \eqref{eq:boundingstabilityconstant}: 
\begin{align}\label{eq:stabilityconstantcounterexamples}
    c_{p,B_R}(f_{\pm})&\leq c\left(1+\frac{2(1+\delta)}{(1-\delta)} \ C_{\mathrm{poinc}}\left(p,B_R,\lvert \mathcal{G} \varphi \rvert^p \right)\right).
\end{align}
By Theorem~B.12 together with Theorem~5.10 in \cite{grohs2019stable}, there exists $\eta$ depending on $p$ but independent of $R>0$ such that 
$$
C_{\mathrm{poinc}}\left(p,B_R,\lvert \mathcal{G} \varphi \rvert^p \right)\leq\eta,
$$
which yields
\begin{align*}
c_{p,B_R}(f_{\pm})&\leq c\left(1+\frac{2(1+\delta)}{(1-\delta)} \eta \right).
\end{align*}
Moreover, by equation \eqref{eq:boundingweights}, we have that 
\[
    \frac{\|\mathcal{G}f_{\pm}\|_{L^p(B_R)}}{\|\mathcal{G}f_{\pm}\|_{L^\infty(B_R)}}
\leq \frac{(1+\delta) \|\mathcal{G}\varphi\|_{L^p(B_R)}}{(1-\delta)\|\mathcal{G}\varphi\|_{L^\infty(B_R)}}
\leq\frac{(1+\delta) \|\mathcal{G}\varphi\|_{L^p(\R^2)}}{(1-\delta)\|\mathcal{G}\varphi\|_{L^\infty(\R^2)}},
\]
as well as
\begin{align*}
    \frac{\|\mathcal{V}_{\varphi'}f_{\pm}\|_{L^\infty(B_R)}}{\|\mathcal{G}f_{\pm}\|_{L^\infty(B_R)}}
    &\leq \frac{\|\mathcal{V}_{\varphi'}\varphi\|_{L^\infty(B_R)}+\gamma\|\mathcal{V}_{\varphi'}\varphi\|_{L^\infty(\R^2)}}{(1-\delta)\|\mathcal{G}\varphi\|_{L^\infty(B_R)}} \\
    &\leq \frac{(1+\delta) \|\mathcal{V}_{\varphi'}\varphi\|_{L^\infty(\R^2)}}{(1-\delta)\|\mathcal{G}\varphi\|_{L^\infty(\R^2)}}.
\end{align*}
Since the constant $c$ in \eqref{eq:stabilityconstantcounterexamples} depends monotonically increasingly on \eqref{eq:max}, the above inequalities allow to upper bound the constant $c$ with a constant $c'$ independent of $R$, $a$ and $\gamma$. Hence, we conclude the proof by defining
\[
    \nu=c'\left(1+\frac{2(1+\delta)}{(1-\delta)} \eta\right), 
\]
which is independent of $R$ and $a$.
\end{proof}

\subsection{Discussion of Theorem~\ref{thm:main}.}

To conclude this section, we give some insights and discuss possible extensions of our main theorem. 
\begin{enumerate}
\item The constant $\nu$ is linked to the stability constant of the Gaussian $\varphi$, which in the result in \cite{grohs2019stable} is estimated by $c(1+\eta)$. The Gaussian $\varphi$ enjoys very strong stability properties for Gabor phase retrieval and the class $M_\nu(B_R)$ with our choice of $\nu$ has stability properties close to that of $\varphi$.

\item It is worth observing that, while Question~\ref{conjecture} is stated for $\Omega=\R^2$, Theorem~\ref{thm:main} is proved for arbitrary large balls $B_R$, with $R>1$. This restriction originates from the bounds on the Poincar\'e constant. However, Theorem~\ref{thm:main} shows that for every sampling rate $a>0$, we can construct functions $f_\pm$ satisfying 
\[
    c_{p,B_R}(f_{\pm})\leq \nu,
\]
with $R>1/a$. The condition $R>1/a$ implies that the ball $B_R$ encloses the centers $(0,0)$ and $(1/a,0)$ of the two bumps of $|\mathcal{G}f_\pm|^p$, and consequently all the features that may affect the local stability constants $c_{p,\R^2}(f_{\pm})$. For this reason, it seems plausible to conjecture that
$
  c_{p,\R^2}(f_{\pm})
$
may also be bounded by a constant $\nu'$ independent of the sampling rate $a$. While a proof of this final argument would allow us to fully answer Question~\ref{conjecture}, this seems to be mainly a technicality.

\item We can extend Theorem~\ref{thm:main} to general lattices of the form $\Lambda=L\Z^2$, $L\in{\mathrm{GL}}_2(\R)$: Given a lattice $\Lambda$, there exist $a>0$ and $\theta\in\R$ such that $\Lambda\subseteq R_{\theta}(\R\times a\Z)$, where $R_{\theta}$ denotes the rotation matrix 
\[
R_{\theta}=\left(\begin{matrix}\cos{\theta} & -\sin{\theta}\\
\sin{\theta} & \cos{\theta}
\end{matrix}\right).
\]
We can therefore adapt the proof of Theorem~\ref{thm:main} to the functions
\[
f^{\theta}_\pm(t)=\mathcal{F}_{-\theta}f_\pm(t),
\]
where $\mathcal{F}_{-\theta}\colon L^2(\R^2)\to L^2(\R^2)$ denotes the fractional Fourier transform\footnote{The \emph{fractional Fourier transform} of order $\theta \in \R$ of a function $f$ is given by
\begin{equation*}
    \mathcal{F}_\theta f (\xi) = \sqrt{1 - \mathrm{i} \cot \theta} \cdot \mathrm{e}^{\pi \mathrm{i} \cot(\theta) \xi^2} \int_\R f(t) \mathrm{e}^{-2\pi\mathrm{i}\left( \csc(\theta) t \xi - \frac{\cot \theta}{2} t^2 \right)} \, \mathrm{d} t, \qquad \xi \in \R,
\end{equation*}
where the square root is defined such that the argument of the result lies in $(-\pi/2,\pi/2]$.} of order $-\theta$. It holds that $f_\pm^{\theta} \in M^{p}(\R)$, $f_+^{\theta}\not\sim f_-^{\theta}$ and
\[
    |\mathcal{G}f_+^{\theta}(x,\omega)|=|\mathcal{G}f_-^{\theta}(x,\omega)|,\quad (x,\omega)\in R_{\theta}(\R\times a\Z)\supseteq\Lambda.
\]
Furthermore, we can directly compute that 
\begin{equation}\label{eq:gaborrotatedandtranslated}
|\mathcal{G}f_\pm^{\theta}(x,\omega)|=|\mathcal{G}f_\pm(R_{-\theta}(x,\omega))|,\quad (x,\omega)\in\R^2.
\end{equation}
The reader may consult \cite{alaifari2021phase,grohs2022foundational} for the detailed proofs of the above facts.
We see from equation~\eqref{eq:gaborrotatedandtranslated} that $|\mathcal{G}f_\pm^{\theta}|$ is the result of a rotation of $|\mathcal{G}f_\pm|$ in the time-frequency plane. Thus, it follows by equation~\eqref{eq:root_set_Gcounterex} that the root sets of $\mathcal{G}f_\pm^{\theta}$ are 
      \begin{equation*}
            \left\{ R_{\theta}\left( \frac{1}{2a} - \frac{a \log \gamma}{\pi}, \pm \frac{a}{2} + 2 a k \right) \,:\, k \in \mathbb{Z} \right\}.
        \end{equation*}
Therefore, given $R > 1$ and $a>0$, the condition 
        \begin{equation*}
            \gamma< \mathrm{e}^{-\frac{\pi}{a}\left( R - \frac{1}{2a} \right)}
        \end{equation*}
ensures that all the roots of $\mathcal{G}f_\pm^{\theta}$ fall outside 
the strip $R_{\theta}((-R,R) \times \R)$ in the time-frequency plane. With this, it is easy to see that analogous arguments as in the proofs of Corollary~\ref{cor:maincor} and  Theorem~\ref{thm:main} apply to $|\mathcal{G}f_\pm^{\theta}|$.

\item A natural question is whether a real-valuedness assumption on the signals combined with a uniform bound on the local Lipschitz constant would restore uniqueness from samples. We expect the answer to this question to be negative: One can construct counterexamples
\[
g_\pm = \varphi\pm\mathrm{i}\gamma \operatorname{M}_{\frac{1}{a}}\varphi\mp\mathrm{i}\gamma \operatorname{M}_{-\frac{1}{a}}\varphi,\quad \gamma>0.
\]
which by \cite[Theorem~3.13]{grohs2022foundational} are real-valued, do not agree up to global phase and satisfy
\[
|\mathcal{G}g_+(x,\omega)|=|\mathcal{G}g_-(x,\omega)|,\quad (x,\omega)\in a\Z\times\R.
\]
Therefore, we see that for any given rectangular lattice $\Lambda$, there exist real-valued functions which are arbitrarily close to the Gaussian and are thus expected to have good local stability properties. At the same time, they do not agree up to global phase but have Gabor transform magnitudes agreeing on $\Lambda$.
\end{enumerate}

\section{Directions of instability}\label{sec:further_remarks}

This section is a general discussion regarding the connection between instability of phase retrieval and Laplacian eigenfunctions and extends some of the ideas in \cite{grohs2019stable}. The main goal of this section is
\begin{enumerate}
    \item to illustrate the idea of Grohs and Rathmair \cite{grohs2019stable} showing that the difficulty of a phase retrieval problem can be captured by a quantity from spectral geometry, the spectral gap of the Laplacian,
    \item then to extend this idea and show that the corresponding eigenfunction indicates where the obstruction in the phase retrieval problem lies,
    \item and to highlight that, by considering the growth of the next few eigenvalues, one can actually show that the ``space of instabilities'' is often finite-dimensional. If $\lambda_k$ is large, then there is at most a $k$-dimensional subspace of functions such that $f$ cannot be stably distinguished from $f + g$ for any $g$ in this subspace.
\end{enumerate}

\subsection{A one-dimensional toy model}
Before discussing the full model in the complex plane, we introduce a very simple one-dimensional toy model for phase retrieval that will illustrate all these ideas using elementary ideas from calculus.

\begin{quote}
    \textbf{Toy Model.} Let $a:[0,1] \rightarrow \mathbb{R}_{>0}$ be a positive, continuous function. Suppose that $f,g:[0,1] \rightarrow \mathbb{R}$ are two smooth real-valued functions such that
    $$ \int_{0}^{1} a(x) (f'(x) - g'(x))^2 dx \qquad \mbox{is small.}$$
    Does this mean that $f(x) \sim g(x) + c$ for some constant $c$?
\end{quote}

It is clear that if $f(x) = g(x) + c$, then the integral is 0. The condition seems to say that $f'(x) \sim g'(x)$ for most points $x$. So, one would expect that, after integrating, the functions only differ by a global unknown constant (playing the role of uniqueness up to the global phase shift in phase retrieval). However, the weight $a(x)$ could become very small in certain regions of space: The derivatives $f'$ and $g'$ could be very different in that small region of space and the difference might hardly be noticeable in the integral.

The analogue of the idea in \cite{grohs2019stable} in this setting is to consider the spectral gap of the second-order differential operator 
\[
    Lu =- \frac{\mathrm{d}}{\mathrm{d} x} \left(a(x) \frac{\mathrm{d}}{\mathrm{d} x} u\right)
\]
which, after integration by parts, can be defined as the largest constant $c > 0$ such that
\[
    \int_0^1 u(x) \,\mathrm{d} x = 0 \implies   \int_0^1 a(x) u'(x)^2 \,\mathrm{d} x \geq c\int_0^1 u(x)^2 \,\mathrm{d} x.
\]
As an easy example, if $a(x) = 1$, then the classical Wirtinger inequality can be phrased as saying
\[
    \int_0^1 u(x) \,\mathrm{d} x = 0 \implies   \int_0^1 u'(x)^2 \,\mathrm{d} x \geq \pi^2\int_0^1 u(x)^2 \,\mathrm{d} x,
\]
so the constant is $c = \pi^2$ in this case.  It is now clear how we can employ this constant in our toy phase retrieval problem: We can normalize both $f$ and $g$ so that they have mean value 0 and then use
\[
    \int_{0}^{1} a(x) (f'(x) - g'(x))^2 \,\mathrm{d} x \geq c \int_0^1 (f(x) - g(x))^2 \,\mathrm{d} x.
\]
If the first integral is small and the constant $c$ is not too small, then we can deduce $f$ and $g$ are close to each other in $L^2$ in this very precise quantitative sense. However, as already mentioned above, problems can arise when $a(x)$ is allowed to be small and this difficulty will be reflected in the spectral gap being small. We will illustrate this with a simple example (see Figure~\ref{fig:new}). In this example, there is a region where $a(x)$ is so small that the integral will not detect big changes in the derivative in that region. 

\begin{center}
    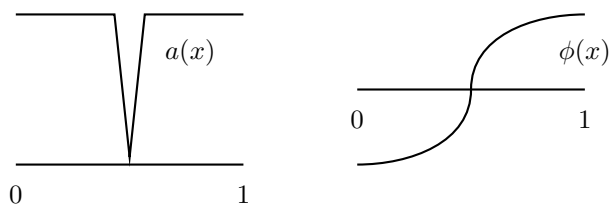
\begin{figure}[h!] \centering
        \begin{tikzpicture}
            \draw [thick] (0,0) -- (3,0);
           \node at (0,-0.4) {0};   \node at (3,-0.4) {1};
           \draw [thick] (0, 2) -- (1.3, 2) -- (1.5, 0.1) -- (1.7, 2) -- (3,2);
           \node at (2.3, 1.5) {$a(x)$};
             \draw [thick] (4.5,1) -- (7.5,1);
              \node at (4.5,1-0.4) {0};   \node at (7.5,1-0.4) {1};
               \draw [thick] (4.5, 0) to[out= 0, in=270] (6,1) to[out=90, in=180] (7.5, 2);
               \node at (7.5, 1.5) {$\phi(x)$};
        \end{tikzpicture}
        \caption{A sketch of an example of a function $a(x)$ (left) corresponding to a difficult ``toy phase retrieval problem''. There is a function $\phi$ (right) for which $\int_0^1 a(x) \phi'(x)^2 \,\mathrm{d} x$ is small even though $\phi$ not close to a constant. }
        \label{fig:new}
    \end{figure}
\end{center}

We can now illustrate our main new idea in this simple case: Further analysis shows that the spectrum of the operator $Lu =- \frac{\mathrm{d}}{\mathrm{d} x} \left(a(x) \frac{\mathrm{d}}{\mathrm{d} x} u\right)$ has one eigenvalue 0 (corresponding to constant functions) and one very small eigenvalue, corresponding to the function $\phi(x)$, \textit{but all the other eigenvalues are actually large}. This means that:
\begin{enumerate}
    \item While the inequality
    \[
        \int_0^1 a(x) u'(x)^2 \,\mathrm{d} x \geq c\int_0^1 u(x)^2 \,\mathrm{d} x
    \]
    is only true for a small constant $c>0$ when looking at all functions with mean value 0, it is actually true with a much larger constant for all functions with mean value 0 that are also orthogonal to $\phi(x)$.
    \item Consequently, while (toy) phase retrieval is difficult for this problem, the only difficulty comes from the function $\phi(x)$. Indeed, the problem is actually easy for all other functions.
    \item There are directions of instability but they are \textit{finite-dimensional} and, in this case, one-dimensional. In particular, if $\int_0^1 a(x) (f'(x) - g'(x))^2 \,\mathrm{d} x$ is small then, in this case, we can deduce in a stable and well-posed manner that $f(x) = g(x) + c_1$ on the left side of the interval and $f(x) = g(x) + c_2$ on the right-hand side. This is the only source of instability.
\end{enumerate}

The implications for phase retrieval are as follows. There are many classical examples of functions where phase retrieval is extremely difficult but they tend to follow the same pattern \cite{alaifari2021Gabor,grohs2019stable}: The signal is strong in two different regions in time-frequency and very weak in between and it is difficult to notice phase shifts in the intermediate region. The implication of our approach is that this is, in a suitable sense, the only difficulty: The phase retrieval problem is actually easy if one is content to recover two different phase shifts (one for each region in time-frequency where the signal is strong), cf. \cite{alaifari2019stable}. We now make this precise.

\subsection{Problem statement and preliminary discussion}

For the purpose of this discussion, we assume that $\Omega \subset \mathbb{C}$ is a bounded domain and that we have two holomorphic functions $F_1, F_2: \Omega \to \mathbb{C}$ where we assume for the sake of simplicity that $\lvert F_1 \rvert > 0$ on all of $\Omega$. The
main question to be discussed is as follows: If
\[
    \lvert F_1 \rvert \approx \lvert F_2 \rvert \mbox{ on most of } \Omega, \mbox{ does this imply that }  F_1 \approx e^{i \alpha} F_2
\]
on most of $\Omega$? Phrased differently: If two holomorphic functions share the same modulus over a large region, does this imply
that one is a global phase-shift of the other? We observe that, throughout this section, the considerations do not invoke the short-time
Fourier transform and are more generally applicable. This subsection may be understood as a short discussion of some of the ingredients in \cite{grohs2019stable} and will set the stage for our subsequent argument. We write
\[
    \inf_{\alpha \in \R} \lVert F_1 - \mathrm{e}^{\mathrm{i} \alpha} F_2 \rVert_{L^2(\Omega)}^2  = \inf_{\alpha \in \R} \int_{\Omega} \left\lvert \frac{F_2(z)}{F_1(z)} - \mathrm{e}^{\mathrm{i} \alpha} \right\rvert^2 \lvert F_1(z) \rvert^2 \,\mathrm{d} z
\]
and thus by defining the measure $\mathrm{d}\mu = \lvert F_1(z) \rvert^2 \,\mathrm{d}z$, we have
\[
    \inf_{\alpha \in \R} \lVert F_1 - \mathrm{e}^{\mathrm{i} \alpha} F_2 \rVert_{L^2(\Omega)}^2
    = \inf_{\alpha \in \R} \int_{\Omega} \left\lvert \frac{F_2(z)}{F_1(z)} - \mathrm{e}^{\mathrm{i} \alpha} \right\rvert^2 \,\mathrm{d} \mu(z).
\]
We can think of the measure $\mu$ as inducing a conformal change of the metric. Assuming $\lvert F_1\rvert$ is sufficiently well behaved, this allows us to interpret the quantity as the $L^2$-norm of a function on a manifold. We now define the real-valued function
\[
    h(z) = \left\lvert \frac{F_2(z)}{F_1(z)} - \mathrm{e}^{\mathrm{i} \alpha} \right\rvert, \qquad z \in \Omega.
\]
At this point, we can invoke the Poincar\'e inequality and argue that
\[
    \int_{\Omega} h(z)^2 \,\mathrm{d}\mu(z) \leq \frac{1}{\mu(\Omega)} \left( \int_\Omega h(z) \,\mathrm{d}\mu(z) \right)^2 + C_{\mathrm{poinc}}(2,\Omega,w)^2 \int_{\Omega} \lvert \nabla h(z) \rvert^2 \,\mathrm{d}\mu(z),
\]
where the Poincar\'e constant, see equation~\eqref{eq:poincdefnsup}, is given by
\[
    C_{\mathrm{poinc}}(2,\Omega,w)^2 = \frac{1}{\lambda_1}
\]
and $\lambda_1$ is the first nontrivial eigenvalue of the Laplace operator on the manifold $(\Omega, \mu)$ equipped with Neumann boundary condition \cite{chavel1984eigenvalues}.
H\"older's inequality immediately implies that
$$  \frac{1}{\mu(\Omega)} \left( \int_\Omega h(z) \,\mathrm{d}\mu(z) \right)^2 \leq  \int_{\Omega} h(z)^2 \,\mathrm{d}\mu(z)$$
with equality if and only if $h$ is constant. In the context of phase retrieval problems considered in this paper, we are mainly interested in the setting where the domain naturally decouples into several subdomains $\Omega = \Omega_1 \cup \Omega_2 \cup \dots \cup \Omega_k$ such that on $\Omega_i$ we have $F_2(z) \sim F_1(z) e^{\alpha_i}$. We observe that in this case, $h$ is not close to a constant globally unless the $\alpha_i$ are all close to each other (corresponding, in essence, to being close to a unified global phase shift).
We also note that if $f : \Omega \to \mathbb{C}$ is analytic in $z_0 \in \Omega$, then
\[
    \left\lvert \nabla \lvert f(z_0)\rvert \right\rvert = \lvert f'(z_0) \rvert
\]
which follows immediately from recalling that the Cauchy--Riemann equations can be geometrically stated as saying that infinitesimal balls are mapped to infinitesimal balls. Therefore, applying this twice, 
\[
    \int_{\Omega} \lvert \nabla h(z) \rvert^2 \,\mathrm{d} \mu(z) = \int_{\Omega}  \left\lvert \left(\frac{F_2(z)}{F_1(z)}\right)' \right\rvert^2 \,\mathrm{d}\mu(z) = \int_{\Omega} \left\lvert \nabla \left\lvert \frac{F_2(z)}{F_1(z)}\right\rvert \right\rvert^2 \,\mathrm{d}\mu(z)
\]
from which we infer
\begin{multline*}
    \inf_{\alpha \in \R} \lVert F_1 - \mathrm{e}^{\mathrm{i} \alpha} F_2 \rVert_{L^2(\Omega)}^2 \leq \frac{1}{\mu(\Omega)} \left( \int_\Omega h(z) \,\mathrm{d}\mu(z) \right)^2 \\
    + C_{\mathrm{poinc}}(2,\Omega,w)^2 \int_{\Omega} \left\lvert \nabla \left\lvert \frac{F_2(z)}{F_1(z)}\right\rvert \right\rvert^2 \,\mathrm{d}\mu(z).
\end{multline*}
If the first term on the right-hand side were to be large, then this would imply that $h$ is typically not small from which we immediately infer that $ F_1 \approx \mathrm{e}^{\mathrm{i} \alpha} F_2$ cannot be true over a large region. So we may henceforth assume that the first term is small.
This leaves us with the second term: If the integral were to be large, then this would be a quantitative measure indicating that $ |F_1| \approx |F_2|$ is not true on most of the domain $\Omega$. However, there is one remaining possibility: It is quite conceivable that the integral is also quite small which would require that $C_{\mathrm{poinc}}(2,\Omega,w)$ is quite large. This in turn implies that $\lambda_1$ is quite small.

\subsection{The Calabi dumbbell example}
\begin{figure}
\centering
\begin{tikzpicture}[scale=0.8]
\draw [very thick] (0,0) to [out=270,in=180] (2,-2) to [out=0,in=180] (3,-0.5) to [out=0,in=180] (4,-2) to [out=0, in=270] (5,0) 
to [out=90, in = 0] (4,2) to [out=180, in = 0] (4,2) to [out=180, in = 0] (3,0.5) to [out=180, in = 0] (2,2) to [out=180, in = 90] (0,0);
\draw[ultra thick, dashed] (3,0.5) to[out=240, in=120] (3,-0.5);
\end{tikzpicture}
\caption{An example of a manifold $(\Omega, \mu)$ isometrically embedded into $\mathbb{R}^2$: This example corresponds to the case where $\lvert F_1(z) \rvert^2$ is large on two separate regions and small everywhere else (including in the area connecting the two regions). This is the prototypical  example of a domain for which $C_{\mathrm{poinc}}(2,\Omega,w)$ is large.}\label{fig:dumbbell}
\end{figure}
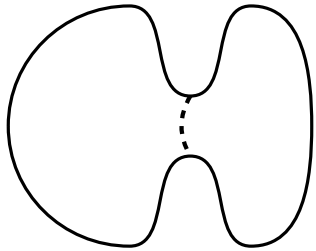

An example is given in Figure \ref{fig:dumbbell}: the classical ``dumbbell'' example is a two-dimensional manifold comprised of two separate regions that are connected via a thin ``bridge''. One way of seeing that the Poincar\'{e} constant for this example is large is to show that $\lambda_1$ is small: Recall that
\[
    \lambda_1 = \inf_{\substack{f \in \mathcal{C}^{\infty}(\Omega) \\ \int_\Omega f \,\mathrm{d}\mu = 0}} \frac{\int_\Omega \lvert \nabla f\rvert^2 \,\mathrm{d}\mu}{\int_\Omega \lvert f \rvert^2 \,\mathrm{d}\mu}.
\]

By taking $f$ to be constant on the left-hand side and right-hand side of the domain and by interpolating linearly in between, we see that $\lvert \nabla f \rvert$ is not necessarily small but that the region over which it is actually nonzero has rather small measure. By making the ``bridge'' thinner, we can make $\lambda_1$ as small as possible. The work of Cheeger then implies that the manifold $(\Omega, \mu)$ can be separated into two distinct parts. Since $\mu = \lvert F_1 \rvert^2$, this simply means that $\lvert F_1 \rvert$ becomes rather small in some regions and this causes the classical and familiar obstruction for phase retrieval: Indeed, when trying to do successful phase retrieval of a function whose information is stored on two separate regions and the function is close to 0 in between, it becomes very difficult to reconstruct the phase globally: Each of the regions may come with a different phase shift. Now, one could wonder whether this is indeed the only obstruction.

\subsection{A refinement}
Our main new idea will be to refine the inequality, valid for all real-valued $h \in W^{1,2}(\Omega, \mu)$,
\[
    \int_{\Omega} h(z)^2 \,\mathrm{d}\mu(z) \leq \frac{1}{\mu(\Omega)} \left( \int_{\Omega} h(z) \,\mathrm{d}\mu(z) \right)^2 + C_{\mathrm{poinc}}(2,\Omega,w)^2 \int_{\Omega} \lvert \nabla h(z) \rvert^2 \,\mathrm{d}\mu(z).
\]
To this end, we introduce a sequence of Laplacian eigenfunctions: These are solutions of
\begin{align*}
    -\Delta u_k &= \lambda_k u_k  \qquad \mbox{inside } (\Omega,\mu),\\
    \frac{\partial u_k}{\partial n} &= 0  \qquad \qquad \mbox{on}~\partial \Omega,
\end{align*}
where
\[
    0 = \lambda_0 < \lambda_1 \leq \lambda_2 \leq \dots
\]
is a discrete sequence of eigenvalues and $n$ is the normal derivative. 

The fundamental inequality that is then being used is
\[
    \int_{\Omega} h(z)^2 \,\mathrm{d}\mu(z) \leq \frac{1}{\mu(\Omega)} \left( \int_{\Omega} h(z) \,\mathrm{d}\mu(z) \right)^2 +  \frac{1}{\lambda_1} \int_{\Omega} \lvert \nabla h(z) \rvert^2 \,\mathrm{d}\mu(z).
\]
We see that in the case of functions with mean value 0, the first term vanishes and we are exactly in the setting described in Section 5.1. The main purpose is now to introduce the following refinement.

\begin{proposition} \label{prop:prop}
    Let $k \in \mathbb{N}$ and let 
    \[
        \pi_k : L^2(\Omega, \mu) \to \operatorname{span}\left\{ u_j \,:\, 1 \leq j \leq k \right\}
    \]
    denote the orthogonal projection. Then, for all real-valued $h \in W^{1,2}(\Omega, \mu)$, we have
    \begin{multline*}
        \int_{\Omega} h(z)^2 \,\mathrm{d}\mu(z) \leq \frac{1}{\mu(\Omega)} \left( \int_{\Omega} h(z) \,\mathrm{d}\mu(z) \right)^2 +  \frac{1}{\lambda_1} \lVert \nabla \pi_k h \rVert_{L^2(\Omega, \mu)}^2 \\
        + \frac{1}{\lambda_{k+1}} \int_{\Omega} \lvert \nabla h(z) \rvert^2 \,\mathrm{d} \mu(z).
    \end{multline*}
\end{proposition}

\begin{proof}
    Let $h \in W^{1,2}(\Omega, \mu)$ be real-valued. We can write 
    \[
        h = \left( h, u_0 \right) u_0 + \sum_{j=1}^{k} \left( h, u_j \right) u_j + \sum_{j = k+1}^\infty \left( h, u_j \right) u_j.
    \]
    We note that this is an orthogonal decomposition into three mutually orthogonal subspaces. The first term gives the constant contribution. As for the other two terms, we note that
    \begin{align*}
        \lVert h \rVert_{L^2(\Omega, \mu)}^2 =  \frac{1}{\mu(\Omega)} \left( \int_{\Omega} h(z) \,\mathrm{d} \mu(z) \right)^2 +  \sum_{j=1}^{k} \left( h, u_j \right)^2  + \sum_{j = k+1}^\infty  \left( h, u_j \right)^2.
    \end{align*}
    By orthogonality, we have
    \[
        \sum_{j=1}^{k} \left( h, u_j \right)^2 = \sum_{j=1}^{k} \left( \pi_k h, u_j \right)^2
    \]
    and thus, arguing as above,
    \[
        \sum_{j=1}^{k} \left( \pi_k h, u_j \right)^2 \leq \frac{1}{\lambda_1} \lVert \nabla \pi_k h \rVert_{L^2(\Omega,\mu)}^2.
    \]
    The same line of reasoning implies 
    \[
        \sum_{j = k+1}^\infty \left( h, u_j \right)^2 \leq \frac{1}{\lambda_{k+1}} \int_{\Omega} \lvert \nabla h(z) \rvert^2 \,\mathrm{d}\mu(z).
    \]
\end{proof}

\subsection{Summarizing the refinement}
Proposition \ref{prop:prop} can now be applied as follows:
\begin{enumerate}
    \item If the spectral gap --- the first nontrivial eigenvalue of the Laplacian $\lambda_1 > 1$ --- is large, then the phase retrieval problem is stable (this is \cite{grohs2019stable}),
    \item and if $\lambda_1$ is small, the problem is generally unstable.
    \item  \textit{However}, if $\lambda_k$ is large, then the phase retrieval problem is actually stable up to a $k$-dimensional subspace (see Proposition \ref{prop:prop}). This means that while phase retrieval is difficult, there are only $k$ different profiles that are difficult to resolve,
    \item and these profiles are given by the Laplacian eigenfunctions.
\end{enumerate}

Moreover, this scenario actually happens in all the typical examples where phase retrieval is difficult: There, usually, $\lambda_1$ is small because of a dumbbell type construction as in Figure~\ref{fig:dumbbell}. However, $\lambda_2$ is already large. The only difficulty in the phase retrieval problem is distinguishing phase shifts on each side of the dumbbell but the problem is otherwise (relatively) easy and one can reconstruct the phase within each region from the magnitude up to one of two possible global phase shifts (one for each side of the dumbbell).

\backmatter

\bmhead{Acknowledgments} The authors would like to thank Alessio Figalli for fruitful discussions on isoperimetric inequalities and Cheeger constants.

\bmhead{Funding} Rima Alaifari, Francesca Bartolucci and Matthias Wellershoff would like to acknowledge funding through the SNSF Grant 200021\_184698. Stefan Steinerberger was partially supported by the NSF (DMS-2123224) and the Alfred P. Sloan Foundation.

\section*{Declaration}

\bmhead{Conflict of interest} 
On behalf of all authors, the corresponding author states that there is no conflict of interest.

\begin{appendices}

\section{List of symbols}\label{app:A}
\begin{itemize}
    \item $\mathcal{V}_\psi f(x,\omega)$ : Short-time Fourier transform defined in equation~\eqref{eqn:STFT}.
    \item $\varphi$ : Normalized Gaussian.
    \item $\mathcal{G} f(x,\omega)$ : Gabor transform defined in equation~\eqref{eqn:Gabor}.
    \item $\mathcal{A}_\Omega(f)$ : Phase retrieval operator defined in Section~\ref{sec:intro}.
    \item $\mathcal{M}_\nu(\mathbb{R}^2)$ : Set of $\nu$-locally stable signals defined in Question~\ref{conjecture}.
    \item $M^p(\mathbb{R})$ : Modulation spaces defined in equation~\eqref{eq:modspace}.
    \item $\lVert f \rVert_{M^p(\mathbb{R})}$ : Norm in the modulation spaces defined in Section~\ref{sec:preliminaries}.
    \item $f \sim g$ : Equivalence up to global phase defined in equation~\eqref{eq:uptoglobalphase}.
    \item $\operatorname{T}_x$ : Translation operators defined in equation~\eqref{eq:definitionoftranslation}.
    \item $\operatorname{M}_\omega$ : Modulation operators defined in Section~\ref{sec:preliminaries}.
    \item $\operatorname{R}_\theta$ : Rotation operators defined in equation~\eqref{eq:definitionofrotation}.
    \item $\mathcal{F}^2(\mathbb{C})$ : Fock space defined in subsection~\ref{ssec:fockspace}.
    \item $\mathcal{B}$ : Bargmann transform defined in subsection~\ref{ssec:fockspace}.
    \item $\Delta$ : Laplace--Beltrami operator.
    \item $0 = \lambda_0 < \lambda_1 < \dots$ : Eigenvalues of the Laplace--Beltrami operator.
    \item $L^p(\Omega,w)$ and $W^{k,p}(\Omega,w)$ : Weighted Lebesgue and Sobolev space defined in subsection~\ref{ssec:prelim_LaplacePoincareCheeger}.
    \item $C_{\mathrm{poinc}}(p,\Omega,w)$: Weighted Poincaré constant defined in equation~\eqref{eq:poincdefnsup}.
    \item $F_\Omega^w$ : Mean of a function $F$ defined in subsection~\ref{ssec:prelim_LaplacePoincareCheeger}.
    \item $w(\Omega)$ : Weight of a set $\Omega$ defined in subsection~\ref{ssec:prelim_LaplacePoincareCheeger}.
    \item $\lVert F \rVert_{\mathcal{D}_{p,q}^{1,4}(\Omega)}$ : Polynomially weighted Sobolev norm defined in equation~\eqref{eq:normmeasuremnts}.
    \item $h_{p,\Omega}(f)$ : Cheeger constant defined in subsection~\ref{ssec:prelim_LaplacePoincareCheeger}.
    \item $\mathfrak{C}(\Lambda)$ : Class of counterexamples defined in Definition~\ref{def:counterexamples}.
    \item $h^\pm(t)$ : Original counterexamples to uniqueness in sampled Gabor phase retrieval defined in equation~\eqref{eq:counterexamples}.
    \item $c_{p,\Omega}(f)$ : Local Lipschitz constant defined in equation~\eqref{eq:local_stability_constant2}.
\end{itemize}

\section{Time-shifting and scaling the counterexamples}
\label{app:time-shift+scale}

In Section~\ref{ssec:counterexamples}, we show that the functions $h_\tau^\pm = \mathcal{B}^{-1} ( \mathrm{e}^{\pi \tau \cdot} \cdot \mathcal{B} h^\pm ) \in L^2(\R)$ --- where $h^\pm$ are defined in equation~\eqref{eq:original_counterexamples} and $\tau \neq 0$ --- are counterexamples to sampled Gabor phase retrieval on $\R \times a \Z$. Additionally, we claim that by applying time-shift and scaling operations to $h_\tau^\pm$, we can obtain the counterexamples $f^\pm_\gamma = \varphi \pm \mathrm{i} \gamma \operatorname{T}_{1/a} \varphi$, where $\gamma > 0$. The proof of this claim is presented in the following.

As a first step, we express $h_\tau^\pm$ as a linear combination of two time-shifted Gaussians. To do so, we recall from Section~\ref{ssec:counterexamples} that the Bargmann transforms of $h^\pm$ are given by 
\begin{equation*}
    \mathcal{B} h^\pm (z) = \mathrm{e}^{\tfrac{\pi}{8 a^2}} \left( \cosh \left(\frac{\pi z}{2a}\right) \pm \mathrm{i} \sinh \left(\frac{\pi z}{2a}\right) \right), \qquad z \in \bbC.
\end{equation*}
Therefore, the relation of the Gabor and the Bargmann transform implies that
\begin{align*}
    \mathcal{G} h^\pm_\tau(x,\omega) &= \mathrm{e}^{-\pi \mathrm{i} x \omega -\frac{\pi}{2}(x^2 + \omega^2)} \mathcal{B} h^\pm_\tau (x- \mathrm{i} \omega) \\
    &= \mathrm{e}^{\pi \tau (x - \mathrm{i} \omega)-\pi \mathrm{i} x \omega -\frac{\pi}{2}(x^2 + \omega^2)} \mathcal{B} h^\pm (x- \mathrm{i} \omega) \\
    &=\begin{multlined}[t]
        \mathrm{e}^{\tfrac{\pi}{8 a^2} + \pi \tau (x - \mathrm{i} \omega) -\pi \mathrm{i} x \omega -\frac{\pi}{2}(x^2 + \omega^2)} \\
       \cdot \left( \cosh \left(\frac{\pi (x-\mathrm{i}\omega)}{2a}\right) \pm \mathrm{i} \sinh \left(\frac{\pi (x-\mathrm{i}\omega)}{2a}\right) \right)
    \end{multlined} \\
    &= \begin{multlined}[t]
        \mathrm{e}^{\tfrac{\pi}{8 a^2} + \pi \tau (x - \mathrm{i} \omega) -\pi \mathrm{i} x \omega -\frac{\pi}{2}(x^2 + \omega^2)} \\
        \cdot \left( \frac{1\mp\mathrm{i}}{2} \mathrm{e}^{-\frac{\pi(x- \mathrm{i} \omega)}{2a}} + \frac{1\pm\mathrm{i}}{2} \mathrm{e}^{\frac{\pi(x - \mathrm{i} \omega)}{2a}} \right)
    \end{multlined} \\
    &=\begin{multlined}[t]
        \frac{1\mp\mathrm{i}}{2} \mathrm{e}^{\frac{\pi}{8a^2} + \frac{\pi}{2}(\frac{1}{2a} - \tau)^2} \mathrm{e}^{-\pi \mathrm{i} (x - \frac{1}{2a} + \tau) \omega} \mathrm{e}^{-\frac{\pi}{2} ((x+\frac{1}{2a}-\tau)^2 + \omega^2)} \\
        + \frac{1\pm\mathrm{i}}{2} \mathrm{e}^{\frac{\pi}{8a^2} + \frac{\pi}{2}(\frac{1}{2a} + \tau)^2} \mathrm{e}^{-\pi \mathrm{i} (x + \frac{1}{2a} + \tau) \omega} \mathrm{e}^{-\frac{\pi}{2} ((x-\frac{1}{2a}-\tau)^2 + \omega^2)},
    \end{multlined}
\end{align*}
for $(x,\omega) \in \R^2$, where we completed the square in the exponent in the final equality.

According to the covariance property of the Gabor transform (cf.~\cite[Lemma~3.1.3 on p.~41]{groechenig2001foundations}) and equation~\eqref{eq:gabortransformgaussian}, it holds that 
\begin{align*}
    \mathcal{G} \operatorname{T}_\alpha \varphi (x,\omega) &= \mathrm{e}^{-2\pi\mathrm{i} \alpha \omega} \mathcal{G} \varphi(x-\alpha,\omega) = \mathrm{e}^{-2\pi\mathrm{i} \alpha \omega} \mathrm{e}^{-\pi\mathrm{i}(x-\alpha)\omega} \mathrm{e}^{-\frac{\pi}{2}((x-\alpha)^2+ \omega^2)} \\
    &= \mathrm{e}^{-\pi\mathrm{i}(x+\alpha)\omega} \mathrm{e}^{-\frac{\pi}{2}((x-\alpha)^2+ \omega^2)},
\end{align*}
for $\alpha \in \mathbb{R}$ and $(x,\omega) \in \mathbb{R}^2$. Therefore, 
\begin{equation*}
    \mathcal{G} h^\pm_\tau = \frac{1\mp\mathrm{i}}{2} \mathrm{e}^{\frac{\pi}{8a^2} + \frac{\pi}{2}(\frac{1}{2a} - \tau)^2} \mathcal{G} \operatorname{T}_{-\frac{1}{2a}+\tau} \varphi + \frac{1\pm\mathrm{i}}{2} \mathrm{e}^{\frac{\pi}{8a^2} + \frac{\pi}{2}(\frac{1}{2a} + \tau)^2} \mathcal{G} \operatorname{T}_{\frac{1}{2a}+\tau} \varphi
\end{equation*}
which implies 
\begin{equation*}
    h^\pm_\tau = \frac{1\mp\mathrm{i}}{2} \mathrm{e}^{\frac{\pi}{8a^2} + \frac{\pi}{2}(\frac{1}{2a} - \tau)^2} \operatorname{T}_{-\frac{1}{2a}+\tau} \varphi + \frac{1\pm\mathrm{i}}{2} \mathrm{e}^{\frac{\pi}{8a^2} + \frac{\pi}{2}(\frac{1}{2a} + \tau)^2} \operatorname{T}_{\frac{1}{2a}+\tau} \varphi,
\end{equation*}
as desired.

Finally, we can time-shift and scale $h^\pm_\tau$ appropriately: In particular, we have
\begin{equation*}
    (1\pm \mathrm{i}) \mathrm{e}^{-\frac{\pi}{8a^2} - \frac{\pi}{2}(\frac{1}{2a} - \tau)^2} \operatorname{T}_{\frac{1}{2a}-\tau} h^\pm_\tau = \varphi \pm \mathrm{i} \mathrm{e}^{\frac{\pi \tau}{a}} \operatorname{T}_{1/a} \varphi
\end{equation*}
which exactly corresponds to $f^\pm_\gamma = \varphi \pm \mathrm{i} \gamma \operatorname{T}_{1/a} \varphi$ with $\gamma = \mathrm{e}^{\pi \tau / a}$. The fact that the functions $f_\gamma^\pm$ are counterexamples to sampled Gabor phase retrieval on $\R \times a \Z$ follows from the invariance of $\mathfrak{C}(\R \times a \Z)$ under time-shifts and scalar multiplication. Alternatively, one may show that $f_\gamma^\pm \in \mathfrak{C}(\R \times a \Z)$ by direct computation (cf.~the proof of Lemma~\ref{lem:new_counterexamples}).

\end{appendices}


\bibliography{sn-bibliography}


\begin{thebibliography}{32}
\ifx \bisbn   \undefined \def \bisbn  #1{ISBN #1}\fi
\ifx \binits  \undefined \def \binits#1{#1}\fi
\ifx \bauthor  \undefined \def \bauthor#1{#1}\fi
\ifx \batitle  \undefined \def \batitle#1{#1}\fi
\ifx \bjtitle  \undefined \def \bjtitle#1{#1}\fi
\ifx \bvolume  \undefined \def \bvolume#1{\textbf{#1}}\fi
\ifx \byear  \undefined \def \byear#1{#1}\fi
\ifx \bissue  \undefined \def \bissue#1{#1}\fi
\ifx \bfpage  \undefined \def \bfpage#1{#1}\fi
\ifx \blpage  \undefined \def \blpage #1{#1}\fi
\ifx \burl  \undefined \def \burl#1{\textsf{#1}}\fi
\ifx \doiurl  \undefined \def \doiurl#1{\url{https://doi.org/#1}}\fi
\ifx \betal  \undefined \def \betal{\textit{et al.}}\fi
\ifx \binstitute  \undefined \def \binstitute#1{#1}\fi
\ifx \binstitutionaled  \undefined \def \binstitutionaled#1{#1}\fi
\ifx \bctitle  \undefined \def \bctitle#1{#1}\fi
\ifx \beditor  \undefined \def \beditor#1{#1}\fi
\ifx \bpublisher  \undefined \def \bpublisher#1{#1}\fi
\ifx \bbtitle  \undefined \def \bbtitle#1{#1}\fi
\ifx \bedition  \undefined \def \bedition#1{#1}\fi
\ifx \bseriesno  \undefined \def \bseriesno#1{#1}\fi
\ifx \blocation  \undefined \def \blocation#1{#1}\fi
\ifx \bsertitle  \undefined \def \bsertitle#1{#1}\fi
\ifx \bsnm \undefined \def \bsnm#1{#1}\fi
\ifx \bsuffix \undefined \def \bsuffix#1{#1}\fi
\ifx \bparticle \undefined \def \bparticle#1{#1}\fi
\ifx \barticle \undefined \def \barticle#1{#1}\fi
\bibcommenthead
\ifx \bconfdate \undefined \def \bconfdate #1{#1}\fi
\ifx \botherref \undefined \def \botherref #1{#1}\fi
\ifx \url \undefined \def \url#1{\textsf{#1}}\fi
\ifx \bchapter \undefined \def \bchapter#1{#1}\fi
\ifx \bbook \undefined \def \bbook#1{#1}\fi
\ifx \bcomment \undefined \def \bcomment#1{#1}\fi
\ifx \oauthor \undefined \def \oauthor#1{#1}\fi
\ifx \citeauthoryear \undefined \def \citeauthoryear#1{#1}\fi
\ifx \endbibitem  \undefined \def \endbibitem {}\fi
\ifx \bconflocation  \undefined \def \bconflocation#1{#1}\fi
\ifx \arxivurl  \undefined \def \arxivurl#1{\textsf{#1}}\fi
\csname PreBibitemsHook\endcsname

\bibitem{prusa2017phase}
\begin{bchapter}
\bauthor{\bsnm{Pr\r{u}\v{s}a}, \binits{Z.}},
\bauthor{\bsnm{Holighaus}, \binits{N.}}:
\bctitle{Phase vocoder done right}.
In: \bbtitle{2017 25th European Signal Processing Conference (EUSIPCO)},
pp. \bfpage{976}--\blpage{980}.
\bpublisher{IEEE},
\blocation{Kos, Greece}
(\byear{2017})
\end{bchapter}
\endbibitem

\bibitem{alaifari2017phase}
\begin{barticle}
\bauthor{\bsnm{Alaifari}, \binits{R.}},
\bauthor{\bsnm{Grohs}, \binits{P.}}:
\batitle{Phase retrieval in the general setting of continuous frames for {B}anach spaces}.
\bjtitle{SIAM Journal on Mathematical Analysis}
\bvolume{49}(\bissue{3}),
\bfpage{1895}--\blpage{1911}
(\byear{2017})
\end{barticle}
\endbibitem

\bibitem{iwen2023phase}
\begin{barticle}
\bauthor{\bsnm{Iwen}, \binits{M.}},
\bauthor{\bsnm{Perlmutter}, \binits{M.}},
\bauthor{\bsnm{Sissouno}, \binits{N.}},
\bauthor{\bsnm{Viswanathan}, \binits{A.}}:
\batitle{Phase retrieval for ${L}^2 ([-\pi, \pi])$ via the provably accurate and noise robust numerical inversion of spectrogram measurements}.
\bjtitle{Journal of Fourier Analysis and Applications}
\bvolume{29}(\bissue{1}),
\bfpage{8}
(\byear{2023})
\end{barticle}
\endbibitem

\bibitem{eldar2014sparse}
\begin{barticle}
\bauthor{\bsnm{Eldar}, \binits{Y.C.}},
\bauthor{\bsnm{Sidorenko}, \binits{P.}},
\bauthor{\bsnm{Mixon}, \binits{D.G.}},
\bauthor{\bsnm{Barel}, \binits{S.}},
\bauthor{\bsnm{Cohen}, \binits{O.}}:
\batitle{Sparse phase retrieval from short-time {F}ourier measurements}.
\bjtitle{IEEE Signal Processing Letters}
\bvolume{22}(\bissue{5}),
\bfpage{638}--\blpage{642}
(\byear{2014})
\end{barticle}
\endbibitem

\bibitem{bojarovska2016phase}
\begin{barticle}
\bauthor{\bsnm{Bojarovska}, \binits{I.}},
\bauthor{\bsnm{Flinth}, \binits{A.}}:
\batitle{Phase retrieval from {G}abor measurements}.
\bjtitle{Journal of Fourier Analysis and Applications}
\bvolume{22}(\bissue{3}),
\bfpage{542}--\blpage{567}
(\byear{2016})
\end{barticle}
\endbibitem

\bibitem{li2017phase}
\begin{barticle}
\bauthor{\bsnm{Li}, \binits{L.}},
\bauthor{\bsnm{Cheng}, \binits{C.}},
\bauthor{\bsnm{Han}, \binits{D.}},
\bauthor{\bsnm{Sun}, \binits{Q.}},
\bauthor{\bsnm{Shi}, \binits{G.}}:
\batitle{Phase retrieval from multiple-window short-time {F}ourier measurements}.
\bjtitle{IEEE Signal Processing Letters}
\bvolume{24}(\bissue{4}),
\bfpage{372}--\blpage{376}
(\byear{2017})
\end{barticle}
\endbibitem

\bibitem{pfander2019robust}
\begin{barticle}
\bauthor{\bsnm{Pfander}, \binits{G.E.}},
\bauthor{\bsnm{Salanevich}, \binits{P.}}:
\batitle{Robust phase retrieval algorithm for time-frequency structured measurements}.
\bjtitle{SIAM journal on imaging sciences}
\bvolume{12}(\bissue{2}),
\bfpage{736}--\blpage{761}
(\byear{2019})
\end{barticle}
\endbibitem

\bibitem{salanevich2023injectivity}
\begin{botherref}
\oauthor{\bsnm{Salanevich}, \binits{P.}}:
Injectivity of multi-window {G}abor phase retrieval.
arXiv:2307.00834
(2023)
\end{botherref}
\endbibitem

\bibitem{alaifari2019stable}
\begin{barticle}
\bauthor{\bsnm{Alaifari}, \binits{R.}},
\bauthor{\bsnm{Daubechies}, \binits{I.}},
\bauthor{\bsnm{Grohs}, \binits{P.}},
\bauthor{\bsnm{Yin}, \binits{R.}}:
\batitle{Stable phase retrieval in infinite dimensions}.
\bjtitle{Foundations of Computational Mathematics}
\bvolume{19},
\bfpage{869}--\blpage{900}
(\byear{2019})
\end{barticle}
\endbibitem

\bibitem{alaifari2021Gabor}
\begin{barticle}
\bauthor{\bsnm{Alaifari}, \binits{R.}},
\bauthor{\bsnm{Grohs}, \binits{P.}}:
\batitle{Gabor phase retrieval is severely ill-posed}.
\bjtitle{Applied and Computational Harmonic Analysis}
\bvolume{50},
\bfpage{401}--\blpage{419}
(\byear{2021})
\end{barticle}
\endbibitem

\bibitem{alaifari2021stability}
\begin{barticle}
\bauthor{\bsnm{Alaifari}, \binits{R.}},
\bauthor{\bsnm{Wellershoff}, \binits{M.}}:
\batitle{Stability estimates for phase retrieval from discrete {G}abor measurements}.
\bjtitle{Journal of Fourier Analysis and Applications}
\bvolume{27},
\bfpage{1}--\blpage{31}
(\byear{2021})
\end{barticle}
\endbibitem

\bibitem{alaifari2021phase}
\begin{barticle}
\bauthor{\bsnm{Alaifari}, \binits{R.}},
\bauthor{\bsnm{Wellershoff}, \binits{M.}}:
\batitle{Phase retrieval from sampled {G}abor transform magnitudes: {c}ounterexamples}.
\bjtitle{Journal of Fourier Analysis and Applications}
\bvolume{28}(\bissue{9}),
\bfpage{1}--\blpage{8}
(\byear{2021})
\end{barticle}
\endbibitem

\bibitem{alaifari2021uniqueness}
\begin{barticle}
\bauthor{\bsnm{Alaifari}, \binits{R.}},
\bauthor{\bsnm{Wellershoff}, \binits{M.}}:
\batitle{Uniqueness of {STFT} phase retrieval for bandlimited functions}.
\bjtitle{Applied and Computational Harmonic Analysis}
\bvolume{50},
\bfpage{34}--\blpage{48}
(\byear{2021})
\end{barticle}
\endbibitem

\bibitem{grohs2023injectivity}
\begin{barticle}
\bauthor{\bsnm{Grohs}, \binits{P.}},
\bauthor{\bsnm{Liehr}, \binits{L.}}:
\batitle{Injectivity of {G}abor phase retrieval from lattice measurements}.
\bjtitle{Applied and Computational Harmonic Analysis}
\bvolume{62},
\bfpage{173}--\blpage{193}
(\byear{2023})
\end{barticle}
\endbibitem

\bibitem{grohs2019stable}
\begin{barticle}
\bauthor{\bsnm{Grohs}, \binits{P.}},
\bauthor{\bsnm{Rathmair}, \binits{M.}}:
\batitle{Stable {G}abor phase retrieval and spectral clustering}.
\bjtitle{Communications on Pure and Applied Mathematics}
\bvolume{72}(\bissue{5}),
\bfpage{981}--\blpage{1043}
(\byear{2019})
\end{barticle}
\endbibitem

\bibitem{grohs2021stable}
\begin{botherref}
\oauthor{\bsnm{Grohs}, \binits{P.}},
\oauthor{\bsnm{Liehr}, \binits{L.}}:
Stable {G}abor phase retrieval in {G}aussian shift-invariant spaces via biorthogonality.
arXiv:2109.02494
(2021)
\end{botherref}
\endbibitem

\bibitem{grohs2022foundational}
\begin{barticle}
\bauthor{\bsnm{Grohs}, \binits{P.}},
\bauthor{\bsnm{Liehr}, \binits{L.}}:
\batitle{On foundational discretization barriers in {STFT} phase retrieval}.
\bjtitle{Journal of Fourier Analysis and Applications}
\bvolume{28}(\bissue{39}),
\bfpage{1}--\blpage{21}
(\byear{2022})
\end{barticle}
\endbibitem

\bibitem{grohs2022non}
\begin{botherref}
\oauthor{\bsnm{Grohs}, \binits{P.}},
\oauthor{\bsnm{Liehr}, \binits{L.}}:
Non-uniqueness theory in sampled {S}{T}{F}{T} phase retrieval.
arXiv:2207.05628
(2022)
\end{botherref}
\endbibitem

\bibitem{grohs2022phaseless}
\begin{botherref}
\oauthor{\bsnm{Grohs}, \binits{P.}},
\oauthor{\bsnm{Liehr}, \binits{L.}}:
Phaseless sampling on square-root lattices.
arXiv:2209.11127
(2022)
\end{botherref}
\endbibitem

\bibitem{grohs2022multi}
\begin{botherref}
\oauthor{\bsnm{Grohs}, \binits{P.}},
\oauthor{\bsnm{Liehr}, \binits{L.}},
\oauthor{\bsnm{Rathmair}, \binits{M.}}:
Multi-window {S}{T}{F}{T} phase retrieval: lattice uniqueness.
arXiv:2207.10620
(2022)
\end{botherref}
\endbibitem

\bibitem{wellershoff2023sampling}
\begin{barticle}
\bauthor{\bsnm{Wellershoff}, \binits{M.}}:
\batitle{Sampling at twice the {N}yquist rate in two frequency bins guarantees uniqueness in gabor phase retrieval}.
\bjtitle{Journal of Fourier Analysis and Applications}
\bvolume{29}(\bissue{1}),
\bfpage{7}
(\byear{2023})
\end{barticle}
\endbibitem

\bibitem{wellershoff2022phase}
\begin{botherref}
\oauthor{\bsnm{Wellershoff}, \binits{M.}}:
Phase retrieval of entire functions and its implications for {G}abor phase retrieval.
arXiv:2202.03733
(2022)
\end{botherref}
\endbibitem

\bibitem{wellershoff2021injectivity}
\begin{botherref}
\oauthor{\bsnm{Wellershoff}, \binits{M.}}:
Injectivity of sampled {G}abor phase retrieval in spaces with general integrability conditions.
arXiv:2112.10136 
(2021)
\end{botherref}
\endbibitem

\bibitem{toft2004continuity}
\begin{barticle}
\bauthor{\bsnm{Toft}, \binits{J.}}:
\batitle{Continuity properties for modulation spaces, with applications to pseudo-differential calculus --- {I}}.
\bjtitle{Journal of Functional Analysis}
\bvolume{207}(\bissue{2}),
\bfpage{399}--\blpage{429}
(\byear{2004})
\end{barticle}
\endbibitem

\bibitem{groechenig2001foundations}
\begin{bbook}
\bauthor{\bsnm{Gr{\"o}chenig}, \binits{K.}}:
\bbtitle{Foundations of Time-frequency Analysis},
\bedition{1}st edn.
\bsertitle{Applied and Numerical Harmonic Analysis}.
\bpublisher{Birkh{\"a}user},
\blocation{Boston, MA}
(\byear{2001})
\end{bbook}
\endbibitem

\bibitem{beneteau2010extremal}
\begin{barticle}
\bauthor{\bsnm{B{\'e}n{\'e}teau}, \binits{C.}},
\bauthor{\bsnm{Carswell}, \binits{B.J.}},
\bauthor{\bsnm{Kouchekian}, \binits{S.}}:
\batitle{Extremal problems in the {F}ock space}.
\bjtitle{Computational Methods and Function Theory}
\bvolume{10},
\bfpage{189}--\blpage{206}
(\byear{2010}).
\bcomment{\url{https://doi.org/10.1007/BF03321762}}
\end{barticle}
\endbibitem

\bibitem{chavel1984eigenvalues}
\begin{bbook}
\bauthor{\bsnm{Chavel}, \binits{I.}}:
\bbtitle{Eigenvalues in {R}iemannian Geometry}.
\bsertitle{Pure and applied mathematics}.
\bpublisher{Academic Press},
\blocation{Orlando, Florida}
(\byear{1984})
\end{bbook}
\endbibitem

\bibitem{cheeger1971lower}
\begin{bbook}
\bauthor{\bsnm{Cheeger}, \binits{J.}}:
\bbtitle{A lower bound for the smallest eigenvalue of the {L}aplacian}.
\bsertitle{Princeton Mathematical Series},
vol. \bseriesno{55},
pp. \bfpage{195}--\blpage{199}.
\bpublisher{Princeton University Press},
\blocation{Princeton, NJ}
(\byear{1971})
\end{bbook}
\endbibitem

\bibitem{titchmarsh1939theory}
\begin{bbook}
\bauthor{\bsnm{Titchmarsh}, \binits{E.C.}}:
\bbtitle{The Theory of Functions},
\bedition{2}nd edn.
\bpublisher{Oxford University Press},
\blocation{Amen House, London E.C.4}
(\byear{1939})
\end{bbook}
\endbibitem

\bibitem{ganapathy1936note}
\begin{barticle}
\bauthor{\bsnm{Iyer}, \binits{V.G.}}:
\batitle{A note on integral functions of order 2 bounded at the lattice points}.
\bjtitle{Journal of the London Mathematical Society}
\bvolume{s1-11}(\bissue{4}),
\bfpage{247}--\blpage{249}
(\byear{1936})
\end{barticle}
\endbibitem

\bibitem{pfluger1937analytic}
\begin{barticle}
\bauthor{\bsnm{Pfluger}, \binits{A.}}:
\batitle{On analytic functions bounded at the lattice points}.
\bjtitle{Proceedings of the London Mathematical Society}
\bvolume{s2-42}(\bissue{1}),
\bfpage{305}--\blpage{315}
(\byear{1937})
\end{barticle}
\endbibitem

\bibitem{dyda2013weighted}
\begin{barticle}
\bauthor{\bsnm{Dyda}, \binits{B.}},
\bauthor{\bsnm{Kassmann}, \binits{M.}}:
\batitle{On weighted {P}oincar{\'e} inequalities}.
\bjtitle{Annales Academiae Scientiarum Fennicae Mathematica}
\bvolume{38},
\bfpage{721}--\blpage{726}
(\byear{2013})
\end{barticle}
\endbibitem

\end{thebibliography}


\end{document}